\documentclass[sts,preprint,noinfoline]{imsart}

\RequirePackage{amsthm,amsmath,amsfonts,amssymb}
\RequirePackage[authoryear]{natbib}
\RequirePackage[colorlinks,citecolor=blue,urlcolor=blue]{hyperref}
\RequirePackage{graphicx}
\RequirePackage{tikz}
\RequirePackage{xcolor}

\makeatletter
\setpkgattr{frontmatter}{style}{\centering}
\setpkgattr{title}{style}{\centering\def\\{\break}}
\setpkgattr{author}{style}{\centering}
\setpkgattr{address}{style}{\centering}
\setpkgattr{title}{skip}{18\p@}
\setpkgattr{authors}{skip}{14\p@}
\setpkgattr{address}{skip}{8\p@}
\setpkgattr{address}{size}{\small\normalfont}

\def\address{\@ifnextchar[{\arxiv@address}{\arxiv@address[]}}
\def\arxiv@address[#1]#2{%
    \smart@par
    \let\author@fmt@init\author@fmt@init@def
    \vskip\address@skip
    \bgroup
        \address@style\address@size\leavevmode
        #2\par
    \egroup
}
\makeatother

\startlocaldefs
\theoremstyle{plain}

\newtheorem{lem}{Lemma}
\newtheorem{cor}{Corollary}
\newtheorem{prp}{Proposition}
\newtheorem{thm}{Theorem}

\theoremstyle{definition}
\newtheorem{exm}{Example}
\newtheorem{dfn}{Definition}

\newtheorem{rmk}{Remark}

\DeclareMathOperator*{\argmax}{arg\,max}

\DeclareMathOperator*{\esssup}{ess\,sup}
\endlocaldefs

\begin{document}

\begin{frontmatter}
\pdfsubject{Preprint}
\title{How e-values generalize hypothesis testing: a Neyman--Pearson lemma for the e-value}
\runtitle{A Neyman--Pearson lemma for the e-value}
\runauthor{Nick W. Koning}

\begin{aug}
\author{\fnms{Nick W.} \snm{Koning}}
\address{%
    Econometric Institute, Erasmus University Rotterdam, the Netherlands\\
    \href{mailto:n.w.koning@ese.eur.nl}{\texttt{n.w.koning@ese.eur.nl}}\\[4pt]
    \today
}
\end{aug}

\begin{abstract}
	While e-values are swiftly rising in prominence, their formal relationship to classical hypothesis testing remains unsettled.
	We develop a unified decision-theoretic framework by viewing the e-value as a multi-decision generalization of a hypothesis test.
	Replacing power by the expected utility of evidence, we derive a Neyman--Pearson lemma for e-values, with classical Neyman--Pearson tests and log-optimal e-values arising from particular utility functions.
	Our main technical contribution is to study such expected utility-optimal e-values without assumptions on the composite null hypothesis and simple alternative.
	As a corollary, we obtain an assumption-free form of the classical Neyman--Pearson lemma for composite null hypotheses, showing power-optimal tests exist and have a likelihood-ratio form.
\end{abstract}

\end{frontmatter}


	\section{Introduction}
		In the traditional Neyman--Pearson framework of hypothesis testing, evidence against a hypothesis is expressed as the decision of a test
		\begin{align*}
			\textnormal{test}: \textnormal{data } \to \{\textnormal{not reject}, \textnormal{reject}\}.
		\end{align*}
		The framework prescribes using a test that decides to falsely reject the hypothesis with probability at most equal to some level of significance $\alpha$:
		\begin{align*}
			\textnormal{Pr}(\textnormal{reject}) \leq \alpha.
		\end{align*}
		A rejection by a test that is valid at a smaller level of significance then indicates stronger evidence, because such a rejection is more surprising if the hypothesis is true.
		
		Compared to the definition of a test, the definition of the e-value may seem arbitrary.
		The e-value is commonly defined as a non-negative random variable $\varepsilon$ that satisfies
		\begin{align*}
			\textnormal{E}[\varepsilon] \leq 1,
		\end{align*}
		if the hypothesis is true \citep{vovk2021values, grunwald2023safe, ramdas2023gametheoretic, ramdas2024hypothesis}.		
		
		We propose a decision-theoretic definition of the e-value that clarifies its relationship to classical testing.
		As a stepping stone towards this definition, we start by incorporating the level $\alpha$ into the decision of the test
		\begin{align*}
			\textnormal{test}: \textnormal{data } \to \{\textnormal{not reject}, \textnormal{reject at level } \alpha \},			
		\end{align*}
		to signify that a rejection at, say, level $1\%$ is a stronger claim than a rejection at level $5\%$.
		We subsequently propose to define an e-value $e$ as a multi-decision generalization of a hypothesis test
		\begin{align*}
			\textnormal{e-value} : \textnormal{data } \to \{\textnormal{not rej.}, \textnormal{rej.  at } \alpha_1, \textnormal{rej.  at } \alpha_2, \dots \},
		\end{align*}
		for various levels $\alpha_1, \alpha_2, ...$
		
		Such a multi-decision generalization of a test requires a generalized notion of validity.
		For this purpose, we rely on recent innovations in extending the Type-I error to data-dependent significance levels \citep{grunwald2022beyond, koning2024post, koning2026extending}.
		Applying this extension to our multi-decision definition of an e-value leads to the choice
		\begin{align*}
			\textnormal{E}[L(\textnormal{e-value})] \leq 1,
		\end{align*}
		for the loss function $L(\textnormal{reject at level } \alpha) = 1 / \alpha$ for every $\alpha$, and $L(\textnormal{not reject}) = 0$.
		
		We recover the common definition of the e-value by identifying the decision to `reject at level $\alpha$' with the numerical value of its loss $L(\textnormal{reject at level } \alpha) = 1 / \alpha$:
		\begin{align*}
			\varepsilon := L(\textnormal{e-value}).
		\end{align*}
		Validity of the e-value then corresponds to $\textnormal{E}[\varepsilon] \leq 1$ if the hypothesis is true.
		A classical test appears as a special $\{0, 1/\alpha\}$-valued e-value.

		\subsection{Continuously interpreting randomized tests}
			Remarkably, we can establish a second connection between e-values and tests.
			A randomized test $\tau$ produces an instruction to reject at level $\alpha$ with a probability in $[0, 1]$ using external randomization.
			
			Suppose we \emph{directly interpret} the instruction to reject with a probability in $[0, 1]$ as evidence against the hypothesis.
			This is at least as informative as the randomized outcome, as one can always choose to apply randomization.
			This means we obtain the benefits of a randomized test, without the downsides of external randomization.
			
			A randomized test $\tau$ is classically said to be valid if its unconditional probability to reject, including randomization, is at most $\alpha$:
			\begin{align*}
				\textnormal{E}[\tau] \leq \alpha.
			\end{align*}
			
			As we no longer intend to randomize, nothing is lost by rescaling by $1/\alpha$ from $[0, 1]$ to $[0, 1/\alpha]$: $\varepsilon_\alpha := \tau / \alpha$.
			This rescaling simultaneously rescales the validity to
			\begin{align*}
				\textnormal{E}[\tau/\alpha]
					= \textnormal{E}[\varepsilon_\alpha]
					\leq 1.
			\end{align*}
			
			This yields a second interpretation of the e-value: \emph{bounded e-values are continuously interpreted randomized tests}.
			Taking $\alpha \to 0$ recovers unbounded e-values.

			The two perspectives on e-values can be united by treating the decision to reject at level $\alpha$ with probability $q$ the same as the decision to reject at level $\alpha / q$:
			\begin{align*}
				\textnormal{reject at level } \alpha \textnormal{ with prob. }q
					\sim \textnormal{reject at level } \alpha/q.
			\end{align*}

		\subsection{Neyman--Pearson lemma for e-values}
			To complete this multi-decision theory of hypothesis testing, we require a generalization of the Neyman--Pearson lemma to multi-decision tests.
			This is necessary, because classical tests maximize the probability to reject at a single level $\alpha$.
			Such a target is inappropriate if we are interested in multi-decision testing.
			
			We generalize the target by specifying our preferences over the different outcomes using a utility function $U$.
			We subsequently maximize the expected utility of the decision under the alternative $Q$:
			\begin{align*}
				\textnormal{E}^{Q}[U(\varepsilon)].
			\end{align*}
			
			Our main technical contribution is to study expected utility-optimal e-values without assumptions on the composite null hypothesis or simple alternative.
			In particular, we provide an expression that verifies an optimizer without assumptions on the utility function $U$.
			This expression is non-decreasing in a likelihood ratio.
			Under mild conditions, we show any optimizer must be of this form.
			Under more regularity conditions, we obtain the explicit form
			\begin{align*}
				\varepsilon^*
					= (U')^{-1}\left(\lambda \frac{dP}{dQ}\right),
			\end{align*}
			where $P$ is an element of an enlargement of the hypothesis $H$ and $\lambda$ is a normalization constant.
			
			Moreover, we provide two flavors of existence results: one that verifies the existence of an optimal e-value \emph{for every pair} $(H, Q)$ under a log-growth condition on $U$, and one that verifies the existence of an optimal e-value \emph{for a particular pair} $(H, Q)$ under a sublinear-growth condition on $U$.
			In addition, we study the uniqueness, exactness and positivity of utility-optimal e-values.
			
			We show that the utility function
			\begin{align*}
				U_\alpha^\textnormal{NP}(x)
					= x \wedge 1/\alpha,
			\end{align*}
			recovers the classical Neyman--Pearson setting, where $x \wedge y = \min(x, y)$.
			This `Neyman--Pearson utility function' does not assign more utility beyond a rejection at some fixed level $\alpha$.
			Applying our results to this setting, we show that a Neyman--Pearson optimal e-value (test) exists without any assumptions on the hypothesis $H$, and that it is of a likelihood ratio test-form.
			
			Log utility $U = \log$ is the near-universal choice in the e-value literature.
			For this setting, we recover the log-optimal or `num\'eraire' e-value results of \citet{larsson2024numeraire} (see also \citet{koolen2022log,grunwald2023safe,lardy2024reverse}).
			We further showcase our results on a capped power utility family, as well as a utility function that specifies relative preferences over common levels such as $1\%$, $5\%$ and $10\%$.
			
		\subsection{Contributions to the e-value literature}
			With this work, we aim to embed the e-value into the classical theory of testing, complementing the relationship between e-values and p-values established in \citet{koning2024post}.
			To the best of our knowledge, we also provide the first unification of Neyman--Pearson optimality theory of tests with the log-optimality theory of e-values \citep{larsson2024numeraire}.		
			
			Our work also contains various critical remarks on the e-value literature.
			Foremost, we aim to end the common view that the e-value is nothing more than a convenient tool to construct binary tests, sequential tests or multiple testing procedures.
			We also challenge the perspective that e-values form a separate paradigm to classical testing.
			Furthermore, we argue that the lauded merging-by-averaging property of e-values is already available in $[0, 1]$-valued testing, though seemingly not used much, which we attribute to an aversion to randomization.
			
			Another goal of our work is to challenge the view that log-utility is a universal panacea, which is how it is often presented in the e-value literature.
			It is so ubiquitous that it is even often referred to as `e-power' \citep{ramdas2024hypothesis}.
			While we see log-utility is an elegant default option, we believe a utility function cannot be universally correct: it must depend on our preferences for evidence within the context of the statistical problem at hand.
			If evidence of at least $1 / \alpha$ is required after one batch of data, then traditional Neyman--Pearson utility $U^{\textnormal{NP}}(x) = x \wedge 1/\alpha$ is an appropriate choice.
			Log-utility is attractive in long-run i.i.d. data collection settings.
			However, these two settings certainly do not exhaust all problems that statisticians encounter.
			While an unfortunate choice of the utility does not invalidate the evidence, there is generally much to gain by tuning the utility function to match the context.
			
			We hope that our work also closes the discussion on how the e-value should be interpreted: an e-value $\varepsilon = e$ should be interpreted as a decision to reject at level $1/e$.
			With the benefit of hindsight, it may be surprising that this is not the common view in the e-value literature.
			For example, in their Chapter 2.7 \citet{ramdas2024hypothesis} cite Jeffreys' rule of thumb to interpret an e-value of 3.16 as similar to a rejection at level $0.05$ instead of a rejection at level $1/3.16 \approx 0.316$.
			We suspect such perspectives stem from a desire to `compensate' for the lossy conversion that happens when e-values are used as a tool to construct a binary level $\alpha$ test: rounding down the e-value to take values in $\{0, 1/\alpha\}$.
			We advise against converting an e-value into a binary test, as this amounts to discarding evidence.
			
		\subsection{Related literature}
			To the best of our knowledge, our Corollary \ref{cor:NP} gives the first assumption-free Neyman--Pearson lemma for an arbitrary composite hypothesis $H$ against a simple alternative $Q$: at every level $\alpha$, a most powerful test exists and has a likelihood-ratio form against an element $P$ of the effective hypothesis $H^\textnormal{eff} \supseteq H$.
			This removes the domination and closure or compactness assumptions of \citet{CvitanicKaratzas2001} and \citet{RudloffKaratzas2010}.
			It also differs from earlier non-dominated extensions, which require either the capacity structure of \citet{HuberStrassen1973}, or derive only a necessary likelihood-ratio form conditional on existence and additional regularity \citep{SunJi2016}.
			The approach to reduce to an element $P$ of the effective hypothesis $H^\textnormal{eff}$ is inspired by \citet{larsson2024numeraire}, following ideas dating back to at least \citet{Lehmann1952} and \citet{Reinhardt1961}.
		
			In mathematical statistics, it is common to study $[0, 1]$-valued tests (see e.g. \citet{van2000asymptotic}).
			Our interpretation is that this literature has already been studying (bounded) e-values, masked by a cosmetic rescaling.
			Compared to this literature, we believe the true innovation of the e-value is an appreciation of the interior $(0, 1)$, which is classically shunned due to its interpretation as an instruction to randomize.
			The appreciation of this interior is mathematically expressed through the utility function $U$: if we stick to Neyman--Pearson utility then the e-value adds nothing over $[0, 1]$-valued testing beyond a different interpretation.
						
			Our work relates to the work of \citet{geyer2005fuzzy}, who also call for a continuous interpretation of a $[0, 1]$-valued test as a `fuzzy' test outcome.
			Besides the fact that they do not rescale to $[0, 1/\alpha]$, the key difference is that they (implicitly) stick to the classical Neyman--Pearson utility.
			As a consequence, their tests are usually not in $(0, 1)$ except in discrete settings.
			We believe this is the reason they present fuzzy tests as an instrument to construct fuzzy confidence sets in discrete data settings.
			
			Surprisingly, we are not aware of prior work that makes this deep connection between e-values and $[0, 1]$-valued tests.
			We suspect that this may be due to the fact that e-values are often contrasted to p-values, and often viewed as a mere `tool' to construct tests.
			To show how the e-value literature struggles with the connection to classical testing we cite the preface of the recent book by \citet{ramdas2024hypothesis}, which first appeared after this work.
			In particular, they define a test as $[0, 1]$-valued yet write: ``[...] a powered test exists if and only if a powered e-value exists'', ``for many [...] problems, the only test we are aware of proceeds by constructing an e-value'', ``anytime-valid inference [...] methods based on e-values and e-processes are the only admissible ones to maintain statistical validity''.
			By arguing that the e-value is a generalization of a test, our work makes the necessity of e-values to any testing problem self-evident.
			
			In subsequent work, \citet{Clerico2024} and \citet{LarssonRamdasRuf2026constraints} study characterizations of e-values for hypotheses generated by constraints.
			Other subsequent work of \citet{LarssonRufRamdas2026} characterizes when a non-trivial test exists for arbitrary composite hypotheses and alternatives, and \citet{RamLarssonRufRamdas2026} derive related results under log-utility.
			This is complementary to our Corollary \ref{cor:NP}, which proves the existence of a test that \emph{attains} optimal power for a composite hypothesis and simple alternative.
						
	\section{From classical tests to e-values}
		\subsection{Background: Hypotheses and tests}
			A hypothesis $H$ is a collection of probability distributions on an underlying sample space $\mathcal{X}$.
			This collection of distributions $H$ can be interpreted as all the data generating processes that share the attribute we intend to falsify.
			
			A test $t$ is classically defined as a function of the underlying data $X$ to the decision space to either reject the hypothesis or not,
			\begin{align*}
				t : \mathcal{X} \to \{\textnormal{do not reject } H, \textnormal{ reject } H\}.
			\end{align*}
			As is common, we suppress the dependence on the underlying data $X$, and write $t$ for $t(X)$, interpreting a test $t$ as a random variable.
			
			Within the classical Neyman--Pearson framework, one should use a test whose probability to reject $H$ is at most equal to some level $\alpha \in (0, 1)$ for every distribution $P$ in the hypothesis:
			\begin{align*}
				P(t = \textnormal{ reject }H) \leq \alpha, \ \textnormal{ for every } P \in H.
			\end{align*}
			
					
		\subsection{The level $\alpha$ in the decision space}
			It is crucial to separate the \emph{nominal} level of a test from its actual probability to falsely reject the hypothesis.
			These are two separate concepts, as we may for example mistakenly claim to have found a rejection at level $1\%$ with a test that has an actual rejection probability of $5\%$.
	
			To capture this distinction, we view the nominal level as a part of the rejection decision:
			\begin{align*}
				t : \mathcal{X} \to \{\textnormal{not reject}, \textnormal{ reject at level } \alpha\}.
			\end{align*}
			This allows us to acknowledge the fact that a rejection at a smaller level is a stronger and distinct claim.
			
			We refer to such a test as a \emph{level $\alpha$ test}, which we say is \emph{valid} if its actual rejection probability is at most $\alpha$:
			\begin{align*}
				P(t = \textnormal{ reject at level } \alpha) \leq \alpha, \ \textnormal{ for every } P \in H.
			\end{align*}
							
		\subsection{Multi-decision tests}
			Distinguishing rejections at different levels enables us to generalize a hypothesis test to a multi-decision test.
			For example, for some pair of levels $\alpha_1, \alpha_2 \in (0, 1)$, we may consider the three-decision test
			\begin{align*}
				e : \mathcal{X} \to \{\textnormal{not reject}, \textnormal{ reject at level } \alpha_1, \textnormal{ reject at level } \alpha_2\}.
			\end{align*}
			More generally, we may consider a multi-decision test that can produce a rejection at any level in $(0, 1)$
			\begin{align*}
				e : \mathcal{X} \to \{\textnormal{not rej.}, \textnormal{ rej. at } \alpha_1, \textnormal{ rej. at } \alpha_2, \dots \}.
			\end{align*}
			
			It remains to define a notion of validity for such a multi-decision generalization of a hypothesis test.
			A naive starting point is to demand that $e$ satisfies
			\begin{align}\label{ineq:validity_naive}
				P(\textnormal{reject at level } \alpha) \leq \alpha,\quad \textnormal{ for every } \alpha \in (0, 1),
			\end{align}
			where `$\textnormal{reject at } \alpha$' means that $e$ makes a decision at least as strong as a rejection at level $\alpha$.
			
			Unfortunately, there is a mismatch between this notion of validity and the nature of a multi-decision test: the level at which a multi-decision test rejects is inherently data-driven, whereas this notion of validity only concerns data-independent levels $\alpha \in (0, 1)$.
			
			To resolve this mismatch, we rely on a recent literature on hypothesis testing with data-dependent levels \citep{grunwald2022beyond, koning2024post}, which studies the validity of a collection of tests $\{t_\alpha\}_{\alpha \in (0, 1)}$ under the post-hoc choice $\widetilde{\alpha}$ of the level (see Remark \ref{rmk:relationship_post-hoc_alpha}).
			To apply these results to our setting, our key insight is that a test at a data-dependent level $t_{\widetilde{\alpha}}$ is itself a multi-decision test, since it may produce a rejection at various levels.
			
			This application leads to a suggested notion of validity that does account for the fact that the decision is inherently data-driven:
			\begin{align}\label{ineq:validity_data_dependent}
				\textnormal{E}^P\left[\frac{P(\textnormal{reject at level } \widetilde{\alpha} \mid \widetilde{\alpha})}{\widetilde{\alpha}}\right] \leq 1, \ \textnormal{ for every } P \in H,
			\end{align}
			for every \emph{data-driven level} $\widetilde{\alpha}$.
			This can be contrasted to \eqref{ineq:validity_naive}, by dividing it by $\alpha$ on both sides:
			\begin{align}\label{ineq:validity_data_independent}
				\frac{P(\textnormal{reject at level } \alpha)}{\alpha} \leq 1, \ \textnormal{ for every } P \in H,
			\end{align}
			for every \emph{data-independent level} $\alpha$.
			Comparing the two shows that \eqref{ineq:validity_data_dependent} requires \eqref{ineq:validity_data_independent} to hold `on average' over the data-driven level.
			
			In Definition \ref{dfn:validity}, we present a more mathematically convenient expression of validity, which is equivalent to \eqref{ineq:validity_data_dependent} (see Theorem 2 in \citet{koning2024post}).
			\begin{dfn}\label{dfn:validity}
				A multi-decision test $e$ is valid for hypothesis $H$ if 
				\begin{align*}
					\textnormal{E}^P\left[L(e)\right] \leq 1, \quad \textnormal{ for every } P \in H,
				\end{align*}
				where $L(\textnormal{rej. at } \alpha) = 1 / \alpha$ and $L(\textnormal{not reject}) = 0$.
			\end{dfn}
			
		\subsection{Connection to the e-value}
			An e-value $\varepsilon$ is usually defined as some non-negative random variable $\varepsilon : \mathcal{X} \to [0, \infty]$ that is valid for the hypothesis $H$ if
			\begin{align*}
				\textnormal{E}^P[\varepsilon] \leq 1, \quad \textnormal{ for every } P \in H.
			\end{align*}
			Unfortunately, this definition can look quite arbitrary on the surface, and we believe it does not do justice to the importance of e-values to hypothesis testing.
			In particular, we argue that \emph{the e-value should be viewed as a multi-decision generalization of a hypothesis test}.
			
			Under this decision-theoretic definition, we recover the standard definition of the e-value as a numerical representation of the multi-decision test $e$:
			\begin{align*}
				\varepsilon := L(e).
			\end{align*}
			This couples the decision to reject at level $\alpha \in (0, 1)$ to the numerical value $1 / \alpha$, and the decision to not reject to the numerical value $0$.
			
			The numerical representation of a classical level $\alpha$ test $t_\alpha$ becomes the $\{0, 1/\alpha\}$-valued binary e-value $\varepsilon_\alpha = L(t_\alpha)$.
			For general multi-decision tests $\varepsilon$, the event $\varepsilon \geq 1/\alpha$ corresponds exactly to the decision to reject at all levels at most $\alpha$.
			
			Note that for $\alpha \in (0, 1)$, we only obtain $\{0\} \cup (1, \infty)$-valued e-values in this manner.
			A $[0, \infty]$-valued e-value can therefore be interpreted as a proper extension of classical tests, adding rejection decisions at `levels $\alpha \geq 1$' that have no interpretation in classical testing.
			A rejection at a level $\alpha > 1$ is not evidence against the hypothesis under the notion of validity in Definition \ref{dfn:validity}, since the decision to always reject at such a level $e \equiv$ `$\textnormal{rej. at } \alpha$' is valid for $\alpha > 1$.
			It must certainly not be mistaken for `evidence against the alternative', without further assumptions.

			\begin{rmk}\label{rmk:relationship_post-hoc_alpha}
				The connection to the e-value that we derive here runs deeper than the connection between post-hoc testing and e-values described in \citet{grunwald2022beyond} and \citet{koning2024post}.
				Indeed, those works connect the e-value to validity of a \emph{collection of binary tests }$\{t_\alpha\}_{\alpha > 0}$ under the post-hoc choice of the level.
				Within our framework, this is analogous to stating that the strongest decision among the tests $\{t_\alpha\}_{\alpha > 0}$ is an e-value.
				However, our framework yields \emph{much} more: even a test at a single data-dependent level $t_{\widetilde{\alpha}}$ is an e-value, since it has a multi-decision outcome: it may yield a rejection at various levels.
			\end{rmk}
	
	\section{From randomized tests to e-values}
		\subsection{Background: randomized tests}\label{sec:randomized}
			In randomized testing, the decision to reject is also allowed to depend on an external source of randomization.
			
			One way to view a randomized test is as a procedure that first uses the data to produce a distribution on the decision space $\{\textnormal{not reject } H, \textnormal{reject } H\}$, and subsequently samples the decision from this distribution.
			As this decision space is binary, the distribution on it is fully captured by the prescribed probability to reject the hypothesis (conditional on the data).
			It is therefore standard in the literature to numerically describe a randomized test as a map that returns this rejection probability
			\begin{align*}
				\tau : \mathcal{X} \to [0, 1].
			\end{align*}
			This conveniently nests non-randomized tests as $\{0, 1\}$-valued, and an outcome $\tau \in (0, 1)$ represents an instruction to reject with probability $\tau$.
			
			Such a randomized test can be operationalized by a binary test $\widetilde{\tau}(X, U)$ that rejects when $U \leq \tau$, where $U \sim \textnormal{Unif}[0, 1]$ independently from the data $X$.
			
			A level $\alpha$ randomized test is said to be valid if its unconditional rejection probability, including both data and randomization, is bounded by $\alpha$:
			\begin{align*}
				(P \times D)(\widetilde{\tau}(X, U) = \textnormal{ reject }H )
					&= (P \times D)(U \leq \tau) \leq \alpha,
			\end{align*}
			for every $P \in H$, where $D = \textnormal{Unif}[0, 1]$.
			Since $D(U \leq \tau) = \tau$, this may be compactly written as
			\begin{align*}
				\textnormal{E}^P[\tau]
					\leq \alpha,\quad\textnormal{ for every } P \in H.
			\end{align*}

			\begin{rmk}[Aversion to randomized testing]
				$[0, 1]$-valued tests are rarely used in practice.
				We attribute this to an aversion to randomization.
				Randomization inherently hampers replicability.
				Moreover, statisticians seem to feel uneasy that the outcome may depend on external randomization instead of only on the data itself.
				
				The `killing blow' is that randomized tests offer almost no advantage in classical testing problems: most powerful tests only require randomization in case of discreteness (see Section \ref{sec:NP_utility}).
				And even then, the probability that randomization is required is often low.
				
				Given the aversion to randomization, it seems that the main reason $[0, 1]$-valued tests are studied is out of mathematical convenience, since this often results in convex optimization problems.
			\end{rmk}

		\subsection{Continuously interpreted randomized tests}
			As a first step towards the e-value, we consider \emph{directly interpreting the outcome of a randomized test $\tau$ in $[0, 1]$ as evidence against the hypothesis.}
			We explicitly propose to not follow this up by external randomization, and so do not arrive at a binary decision.
			This is equivalent to directly interpreting the randomization distribution on $\{\textnormal{not reject }H, \textnormal{reject }H\}$ as our decision.
			
			As this continuous interpretation does not involve randomization, we refer to these tests as $[0, 1]$-valued tests instead of randomized tests.
			We stress that a $[0, 1]$-valued test generalizes a binary test, as the binary decision is a special case $\{0, 1\} \subset [0, 1]$.
			Equivalently, binary tests can be viewed as mapping the data to a distribution that assigns a point mass of 1 on one of the options in $\{\textnormal{not reject }H, \textnormal{reject }H\}$.
			
			We couple the interpretation of a $[0, 1]$-valued test to our intuitions about randomized testing: it expresses the willingness of the test (given the data) to reject the hypothesis with a certain probability.
			This also reveals that a $[0, 1]$-valued test is at least as informative as a random decision: if desired, we may always choose to apply external randomization to retrieve a binary decision.
			Moreover, if a $[0, 1]$-valued test equals 1, then we can directly interpret it as a rejection of the hypothesis.
			
			As with a randomized test, we say that a $[0, 1]$-valued test is valid if
			\begin{align*}
				\textnormal{E}^P[\tau]
					\leq \alpha,\quad\textnormal{ for every } P \in H.
			\end{align*}
			
			\begin{rmk}[Tests as conditional significance levels]
				In subsequent work, \citet{koning2025sequentializing} provide an additional interpretation of a $[0, 1]$-valued test.
				In particular, after observing the outcome of a valid level $\alpha$ test $\tau_1$, we may initialize a second test $\tau_2$ that is valid \emph{at conditional significance level } $\tau_1$.
				That is, $\textnormal{E}^P[\tau_2 \mid \tau_1] \leq \tau_1$, $P$-a.s. for every $P \in H$.
				The overall procedure then remains unconditionally valid at the original level $\alpha$, since $\textnormal{E}^P[\tau_2] = \textnormal{E}^P[\textnormal{E}^P[\tau_2 \mid \tau_1]] \leq \textnormal{E}^P[\tau_1] \leq \alpha$, for every $P \in H$.
				
				A consequence of this perspective is that a test outcome in $[0, 1]$ can be interpreted as the `current significance level' at which we are allowed to do a test while maintaining validity at the unconditional level $\alpha$.
				This shows that such outcomes in $(0, 1)$ have evidential value, since they can be carried forth to subsequent testing.
				
				Using a test outcome as a conditional significance level can also be interpreted as a generalization of randomized testing, where randomized testing can be interpreted as conducting the second test $\tau_2$ with uninformative data, so that the rejection probability is equal to the conditional level $\tau_1$.
			\end{rmk}

		\subsection{Connecting to the e-value by rescaling}
			The choice to numerically represent a randomized test as $[0, 1]$-valued was arbitrary, primarily based on its convenient interpretation as a rejection probability.
			As we no longer intend to randomize, there is little harm in rescaling.
			In particular, without any loss of generality, we can rescale by $1/\alpha$, instead defining a level $\alpha$ test as
			\begin{align*}
				\varepsilon_\alpha = \tau / \alpha,
			\end{align*}
			changing its codomain from $[0, 1]$ to $[0, 1/\alpha]$.
			A benefit of this rescaling is that it clearly differentiates the strength of rejection decisions across different levels: rejections at level 5\% and 1\% now appear numerically as $20$ and $100$, instead of being both equal to $1$.
			
			With the rescaling, we must also rescale the validity of our now $[0, 1/\alpha]$-valued numerical representation of our test $\varepsilon_\alpha$ to
			\begin{align*}
				\textnormal{E}^P[\varepsilon_\alpha]
					\leq 1, \quad \textnormal{ for every } P \in H.
			\end{align*}
			
			Applying this cosmetic rescaling, we arrive at the common definition of the e-value: bounded e-values \emph{are} continuously-interpreted level $\alpha$ randomized tests, thinly veiled by rescaling.
			Moreover, the class of all level $\alpha \in (0, 1]$ randomized tests is equivalent to the class of all bounded e-values.
			
			To go beyond bounded e-values, we can repeat the same operations taking $\alpha \to 0$, lifting the domain of the e-value to $[0, \infty)$ or even $[0, \infty]$.
			In the $\alpha = 0$ limit we lose the interpretation that $\alpha \varepsilon_\alpha$ is the rejection probability of a randomized test (see Remark \ref{rmk:alpha_to_0}).
			Overall, this means that the only benefit of the e-value over $[0, 1]$-valued tests is lifting the upper bound, which seems to be of little practical importance.
			
			\begin{rmk}[Level 0 test on original scale]\label{rmk:alpha_to_0}
				Taking the limit as $\alpha \to 0$ yields something richer than what one would classically call a level $0$ test on the $[0, 1]$-scale.
				Indeed, there a level $0$ test $\tau_0$ would classically be said to be valid if $\sup_{P \in H} \textnormal{E}^P[\tau_0] \leq 0$.
				This is equivalent to $\tau_0 = 0$, $P$-almost surely for every $P \in H$.
				This means $\tau_0$ can only be positive on a set of zero probability for every $P \in H$.
				Moreover, on such an `$H$-null' set, $\tau_0$ is unrestricted so it is reasonable to set it to 1 on such a set (assuming we prefer rejections if the null is false).
				This characterizes level 0 tests on the $[0, 1]$-scale: they are $0$ except for some $H$-null set, on which they are arbitrary.
				This would correspond to a $\{0, \infty\}$-valued e-value, while taking the limit $\alpha \to 0$ yields the richer $[0, \infty]$-valued e-value.
			\end{rmk}

			\begin{rmk}[Betting and rescaling]\label{rmk:betting_rescaling}
				The rescaling has a nice betting interpretation in the spirit of \citet{shafer2021testing}.
				In particular, the traditional $[0, 1]$-scale can be interpreted as betting towards some desired target wealth normalized to $1$, starting from a capital $\alpha$ as a fraction of the desired wealth.
				This means that rejecting at level $\alpha$ is akin to reaching a wealth of $1$ starting from $\alpha$.
				The rescaling we discuss in this section can be interpreted as instead normalizing the starting capital to 1.
				This makes setting a desired target wealth optional.
			\end{rmk}
			
	\section{Converting, merging and sequential}
		\subsection{Converting e-values into tests}
			It is common to convert an e-value into a non-random or randomized test, using the conversions $\varepsilon \mapsto \textnormal{I}\{\varepsilon \geq 1/\alpha\}$ and $\varepsilon \mapsto \alpha\varepsilon \wedge 1$.
			The latter can be operationalized as $\varepsilon \mapsto \textnormal{I}\{\varepsilon \geq U/\alpha\}$, with $U \sim \textnormal{Unif}[0, 1]$.
			The validity of such tests is then motivated through Markov's inequality or the recently introduced randomized Markov's inequality \citep{ramdas2023randomized}.
			
			The connection we establish between e-values and tests bypasses Markov's inequality.
			The conversion of an e-value into a binary test can be interpreted as simply rounding down the e-value to the domain $\{0, 1/\alpha\}$ or $[0, 1/\alpha]$.
			This is akin to simply discarding any evidence above $1/\alpha$ in the randomized setting, and additionally discarding all evidence in $(0, 1/\alpha)$ in the binary setting.
			In the betting interpretation, this is equivalent to throwing away excess wealth, which seems ill-advised.
			
			Such discarding of evidence only seems defensible if one truly does not desire anything beyond a rejection at the specified level $\alpha$.
			However, in that case one would have likely been better off using an e-value that was confined to $\{0, 1/\alpha\}$ in the first place: a level $\alpha$ test.
			
			Overall, we recommend directly reporting the e-value and avoiding discarding evidence as much as possible.
			
			\begin{prp}\label{prp:e-value_to_test}
				Let $\alpha \in (0, 1)$.
				If $\varepsilon$ is an e-value, then
				\begin{itemize}
					\item $\varepsilon \wedge 1/\alpha$ is a level $\alpha$ randomized test,
					\item $\textnormal{I}\{\varepsilon \geq 1/\alpha\}/\alpha$ is a level $\alpha$ test.
				\end{itemize}
				If $\varepsilon$ is valid for a hypothesis $H$, then both are valid for $H$.
			\end{prp}
			\begin{proof}
				Note that $\textnormal{I}\{\varepsilon \geq 1/\alpha\}/\alpha = \lfloor \alpha \varepsilon \wedge 1 \rfloor / \alpha$, so that both operations simply weaken the rejection claims down to the relevant domains.
				Moreover,  $\varepsilon \wedge 1/\alpha \leq \varepsilon$ and $\lfloor \alpha \varepsilon \wedge 1 \rfloor / \alpha \leq \varepsilon$, so that validity of $\varepsilon$ implies validity of the tests.
			\end{proof}
			
			\begin{rmk}[Deterministic Markov's Inequalities]
				The proof of Proposition \ref{prp:e-value_to_test} is deeply related to the `Deterministic Markov's Inequalities' introduced in \citet{koning2024post}.
				For a non-negative value $X$, and an independent $U \sim \textnormal{Unif}[0, 1]$, they are given by:
				\begin{align*}
					\textnormal{I}\{\alpha X \geq 1\}
						&= \lfloor \alpha X \wedge 1 \rfloor
						\leq \alpha X \wedge 1 \\
						&= \mathbb{P}_U(\alpha X \geq U)
						\leq \alpha X.
				\end{align*}
				While these inequalities may seem trivial, applying the expectation over $X$ to all terms recovers Markov's Inequality $\mathbb{P}_X[\alpha X \geq 1] \leq \alpha \mathbb{E}_X[X]$, the Randomized Markov's Inequality of \citet{ramdas2023randomized} $\mathbb{P}_{U, X}[\alpha X \geq U] \leq \alpha\mathbb{E}_{X}[X]$, and a strengthening of Markov's Inequality $\mathbb{P}_X[\alpha X \geq 1] \leq \mathbb{E}_{X}[\alpha X \wedge 1]$ that works for possibly non-integrable non-negative random variables $X$.
				The latter is related to the work of \citet{wang2024extended} who implicitly use it to derive a version of Ville's Inequality for non-integrable non-negative (super)martingales and $e$-processes.
				
				Lastly, a key idea in \citet{koning2024post} is to observe that the inequalities can be made to hold with equality by dividing by $\alpha$, and then taking the supremum over $\alpha$:
				\begin{align*}
					\sup_{\alpha > 0} \frac{\textnormal{I}\{\alpha X \geq 1\}}{\alpha}
						= \sup_{\alpha > 0} X \wedge 1/\alpha
						= X.
				\end{align*}
				Taking the expectation over $X$ yields a `Markov's Equality'.
			\end{rmk}

		\subsection{Merging}
			An advantage of e-values is that they are closed under averaging, unlike classical binary tests.
			Indeed, suppose we have two binary level $\alpha$ tests, numerically represented as having $\{0, 1/\alpha\}$-valued outcomes.
			Then their average is a multi-decision test taking value in $\{0, .5 / \alpha, 1 / \alpha\}$.
			This means that classical binary tests are not closed under such averaging.
			
			In contrast, any (weighted) average of valid $[0, \infty]$-valued multi-decision tests remains a valid multi-decision test, since it is still $[0, \infty]$-valued and
			\begin{align*}
				\textnormal{E}^P[w \varepsilon_1 + (1 - w) \varepsilon_2]
					&= w\textnormal{E}^P[\varepsilon_1] + (1 - w) \textnormal{E}^P[\varepsilon_2] \\
					&\leq w + (1-w)
					= 1,
			\end{align*}
			for every $P \in H$, $w \in [0, 1]$, if $\varepsilon_1$ and $\varepsilon_2$ are valid for $H$.
			
			Such a closure under averaging also holds for $[0, 1]$-valued tests: the average of two such tests that are valid at level $\alpha_1$ and $\alpha_2$ is valid at the average of their levels $(\alpha_1 + \alpha_2) / 2$.
			While possible, this merging does not seem to be widely used in practice.
			We believe this may be due to the fact that this more often yields $(0, 1)$-valued outcomes, which is disliked due to its interpretation as an instruction to randomize.

		\subsection{Sequential testing}
			E-values are frequently presented in the literature as `enabling' sequential testing.
			We use our new insights to argue that the e-value does not enable sequential testing, but yields a multi-decision form of sequential testing.
			
			To describe sequential testing, we require a filtration $(\mathcal{F}_t)_{t \geq 0}$, $\mathcal{F}_{t + 1} \supseteq \mathcal{F}_t$, where $\mathcal{F}_t$ describes the available information at time $t \geq 0$.
			We can then define a level $\alpha$ test process as a sequence of tests $(\varepsilon_\alpha^t)_{t \geq 0}$, where each test $\varepsilon_\alpha^t$ is $\{0, 1/\alpha\}$-valued and $\mathcal{F}_t$-measurable: known at time $t$.
			
			The right notion of validity for such a test process is \emph{anytime validity}:
			\begin{align*}
				\varepsilon_\alpha^{\widetilde{t}} \textnormal{ is valid for every \emph{stopping time} } \widetilde{t}.
			\end{align*}
			Upgrading from tests to e-values, we obtain a multi-decision generalization of a test process $(\varepsilon^t)_{t \geq 0}$, where every $\varepsilon^t$ is $[0, \infty]$-valued.
			Such a multi-decision test process is often referred to as an \emph{e-process}, and similarly said to be anytime valid if
			\begin{align*}
				\varepsilon^{\widetilde{t}} \textnormal{ is valid for every \emph{stopping time} } \widetilde{t}.
			\end{align*}
			
			An e-process offers a great advantage over a classical test process: a (non-dominated) anytime valid level $\alpha$ test process is necessarily of the binary form $0$, $0$, $\dots$, $0$, $1/\alpha$, $1/\alpha$, $\dots$, where the switch from a non-rejection to a rejection at level $\alpha$ happens at some stopping time $\widetilde{t}$.
			In contrast, an e-process exhibits how evidence against the hypothesis continuously progresses over time.
			
			We believe that the e-value is commonly presented as `enabling' sequential testing, because it is often easier to first construct an e-process and then round it down to $\{0, 1/\alpha\}$ by discarding excess evidence:
			\begin{align*}
				\varepsilon^t \mapsto \textnormal{I}\{\varepsilon^t \geq 1/\alpha\} / \alpha,
			\end{align*}
			which is anytime valid as $\textnormal{I}\{\varepsilon^t \geq 1/\alpha\} / \alpha \leq \varepsilon^t$.
			As stated before, we generally recommend against such discarding of evidence.
		
	\section{A Neyman--Pearson lemma for the e-value}
		We now shift our focus to what happens if the hypothesis $H$ is false: if some alternative distribution $Q$ is true.
	
		For a binary level $\alpha$ test, the distribution of its outcome under the alternative $Q$ is completely characterized by its rejection probability.
		It is therefore natural to attempt to maximize this rejection probability under the alternative.
		
		An e-value may return a rejection at various different levels, so that there is no natural power objective: one must specify preferences over different amounts of evidence.
		
		We express these preferences through a measurable utility function
		\begin{align*}
			U : [0, \infty] \to [-\infty, \infty].
		\end{align*}
		We then consider maximizing the expected utility of evidence under the alternative $Q$:
		\begin{align}\label{obj:expected_utility_naive}
			\textnormal{E}^Q[U(\varepsilon)].
		\end{align}
		Given an alternative $Q$, the choice of the utility function $U$ can be interpreted as a risk-reward trade-off: do we wish to reliably gather some evidence, or gamble for a lot of evidence with low probability?
		

		\subsection{Utility optimality}
			In this section, we study the properties of expected utility-optimal e-values against $Q$ for a simple hypothesis $H = \{P\}$, containing a single probability measure $P$.
			
			We write 
			\begin{align*}
				Q = Q_a + Q_s, \quad Q_a \ll P, \quad Q_s \perp P,
			\end{align*}
			for the Lebesgue decomposition of $Q$ with respect to $P$.
			For readers unfamiliar with the Lebesgue decomposition, $Q_a$ can (roughly) be interpreted as the restriction of $Q$ to the region where $P$ is positive, and $Q_s$ where $P$ is zero.\footnote{More precisely, there exists a measurable set $N$ such that $P(N) = 0$, $Q_s(N^c) = 0$ and $Q_a(N) = 0$.}
			
			The expected-utility objective \eqref{obj:expected_utility_naive} requires some refinement.
			In particular, the behavior of an e-value on the singular part $Q_s$ can be chosen freely without affecting its validity under $P$. 
			We therefore only consider the meaningful comparison under the absolutely continuous part $Q_a$.
			
			\begin{dfn}\label{dfn:U-optimal}
				A valid e-value $\varepsilon^*$ for $H$ is said to be $U$-optimal against $Q_a$ if, for every valid e-value $\varepsilon$ for $H$,
				\begin{align}\label{ineq:U-optimal}
					\textnormal{E}^{Q_a}[U(\varepsilon^*) - U(\varepsilon)] \geq 0,
				\end{align}
				and the left-hand side is well-defined.
			\end{dfn}
			With `well-defined', we mean that the integrand avoids an ambiguous `$\infty - \infty$' $Q_a$-almost surely, and that either its positive or negative part has finite expectation.
			
			Replacing $Q$ by $Q_a$ avoids that the objective is blind to what happens under $Q_a$ if $U(\varepsilon^*) = \infty$ with positive $Q_s$-probability.	
			Throughout, one can interpret $\varepsilon^*$ to be chosen such that $\varepsilon^*(x) \in \argmax_{z \in [0, \infty]} U(z)$, $Q_s$-almost surely, assuming a maximum is attained.
			Moreover, we put both utility functions under the same expectation to enable some comparisons where the two expected utilities are not separately well-defined.
			
		\subsection{Neyman--Pearson lemma for e-values}
			In Theorem \ref{thm:EU_NP}, we present our `Neyman--Pearson lemma for the e-value', for an arbitrary utility function $U$.
			A proof may be found in the appendix, alongside all other proofs.
			In the result, we use the convention $0 \times U(z) = 0$.
						
			\begin{thm}\label{thm:EU_NP}
				Let $U : [0, \infty] \to [-\infty, \infty]$ be measurable.
				Suppose that for some $\lambda \in [0, \infty)$, $\varepsilon^*$ is a valid e-value that satisfies
				\begin{align}\label{opt:U-optimal_e-value}
					\varepsilon^*(x)
						\in \argmax_{z \in [0, \infty)} \frac{dQ_a}{dP}(x) U(z) - \lambda z,\quad P\textnormal{-a.s.},
				\end{align}	
				with $\lambda (1 - \textnormal{E}^P[\varepsilon^*]) = 0$, and $U(\varepsilon^*)$ finite $Q_a$-almost surely.
				Then $\varepsilon^*$ is $U$-optimal.
			\end{thm}
			
			Next, we present a converse to Theorem \ref{thm:EU_NP}.
			This converse shows that the pointwise sufficient condition in Theorem \ref{thm:EU_NP} is also necessary for $U$-optimality, under mild conditions on $U$.
		
			\begin{prp}[Converse to Theorem \ref{thm:EU_NP}]\label{prp:EU_converse}
				Assume that $U$ is concave and non-decreasing, and finite on $(0, \infty)$.
				Suppose $\varepsilon^*$ is $U$-optimal.
				Then there exists a $\lambda \in [0, \infty)$ such that $\varepsilon^*$ satisfies \eqref{opt:U-optimal_e-value} and $\lambda (1 - \textnormal{E}^P[\varepsilon^*]) = 0$.
			\end{prp}
			
			The main gap between Proposition \ref{prp:EU_converse} and Theorem \ref{thm:EU_NP} is concavity, since non-decreasingness and finiteness on $(0, \infty)$ seem fairly natural in this context.
			A way to `close' this gap is to replace $U$ by its concave majorant $\overline{U}$.
			As shown in Proposition \ref{prp:concavification}, any optimizer in Theorem \ref{thm:EU_NP} remains an optimizer in this concave problem, but concavification may add optimizers.
			
			\begin{prp}[Concavification]\label{prp:concavification}
				For $a,b \in [0, \infty)$,
				\begin{align*}
					\sup_{z \in [0,\infty)}\left\{aU(z) - bz\right\}
						= \sup_{z \in [0,\infty)}\left\{a\overline{U}(z) - bz\right\}.
				\end{align*}
				Moreover, every maximizer of the first problem is also a maximizer of the second.
			\end{prp}
			
			\begin{rmk}[Non-decreasing in LR]
				A key consequence of \eqref{opt:U-optimal_e-value} is a likelihood ratio-ordering: if $\lambda > 0$, then a larger value of $dQ_a/dP$ cannot correspond to a smaller value of $\varepsilon^*$.
				If the pointwise maximizer is unique, $\varepsilon^*$ is therefore a non-decreasing function of $dQ_a/dP$.
				By Proposition \ref{prp:EU_converse}, this ordering also applies to every $U$-optimal e-value under mild assumptions.
			\end{rmk}
			
			\begin{rmk}[$Q_a$-a.s. finite assumption]
				The assumption in Theorem \ref{thm:EU_NP} that $U(\varepsilon^*)$ is finite $Q_a$-almost surely ensures that the relative expected utility in \eqref{ineq:U-optimal} is well-defined.
				Lemma \ref{lem:well-defined} shows that under the pointwise optimality condition \eqref{opt:U-optimal_e-value}, this is in fact equivalent to the well-definedness requirement.
				We give a sufficient condition in Lemma \ref{lem:well-defined_U}.
			\end{rmk}
			
			\begin{lem}[Well-defined objective]\label{lem:well-defined}
				Let $\varepsilon^*$ be a valid e-value.
				If $\textnormal{E}^{Q_a}[U(\varepsilon^*) - U(\varepsilon)]$ is well-defined for every valid e-value $\varepsilon$, then $U(\varepsilon^*)$ is finite $Q_a$-almost surely.
				
				Conversely, suppose that for some $\lambda \in [0,\infty)$, $\varepsilon^*$ satisfies \eqref{opt:U-optimal_e-value}.
				If $U(\varepsilon^*)$ is finite $Q_a$-almost surely, then $\textnormal{E}^{Q_a}[U(\varepsilon^*) - U(\varepsilon)]$ is well-defined for every valid e-value $\varepsilon$.
			\end{lem}
			
			\begin{lem}[Conditions on $U$ for well-definedness]\label{lem:well-defined_U}
				Suppose that $\lambda \in [0,\infty)$ and a valid e-value $\varepsilon^*$ together satisfy \eqref{opt:U-optimal_e-value}.
				Assume that $U(z) < \infty$ for every finite $z$, and $U(z) > -\infty$ for at least one finite $z$.
				Then $\textnormal{E}^{Q_a}[U(\varepsilon^*) - U(\varepsilon)]$ is well-defined for every valid e-value $\varepsilon$.
			\end{lem}
			
		\subsection{Existence}
			While Theorem \ref{thm:EU_NP} and Proposition \ref{prp:EU_converse} characterize utility-optimal e-values, they do not guarantee their existence.
			In this section, we study the existence both for a particular pair $(P, Q)$ and uniformly in all pairs $(P, Q)$.
			
			In Theorem \ref{thm:existence_simple}, a sublinear growth condition on $U$ ensures that the pointwise optimization problem \eqref{opt:U-optimal_e-value} has a measurable solution.
			For a particular pair $(P,Q)$, existence then follows if these pointwise solutions are $P$-integrable.
			
			Corollary \ref{cor:existence_log_growth} subsequently considers a stronger logarithmic growth condition on $U$ that guarantees this integrability, and thereby existence, for \emph{every} pair $(P, Q)$.
			If $U$ is differentiable, a sufficient condition for \eqref{ineq:log_growth} is $\sup_{z \geq 0} z U'(z) < \infty$.
			
			\begin{thm}[Pair-specific existence]\label{thm:existence_simple}
				Suppose that $U$ is concave, upper semicontinuous and finite on $(0, \infty)$.
				Assume that
				\begin{align}\label{ineq:sublinear}
					\limsup_{z \to \infty} \frac{U(z)}{z}
						\leq 0.
				\end{align}
				Then, for every $\lambda > 0$, \eqref{opt:U-optimal_e-value} has a non-empty compact interval of maximizers.
				Its smallest maximizer $\varepsilon_\lambda$ is measurable and non-increasing in $\lambda$.
				
				If $\textnormal{E}^P[\varepsilon_\lambda] < \infty$ for $\lambda > 0$, then there exists a pair $(\lambda^*, \varepsilon^*)$ satisfying the conditions of Theorem \ref{thm:EU_NP}, so that a $U$-optimal e-value exists.
			\end{thm}
			
			\begin{cor}[Existence under logarithmic growth]\label{cor:existence_log_growth}
				Suppose that $U$ is concave, upper semicontinuous and finite on $(0, \infty)$, and assume there exists some $K \geq 0$ such that 
				\begin{align}\label{ineq:log_growth}
					U(y) - U(x) \leq K \log(y/x),
				\end{align}
				for $1 \leq x \leq y < \infty$.
				Then there exists a pair $(\lambda^*, \varepsilon^*)$ satisfying the conditions of Theorem \ref{thm:EU_NP}, so that a $U$-optimal e-value exists.
			\end{cor}

		\subsection{Structure under additional assumptions on $U$}
			In this section, we describe properties of optimal e-values that are implied by the structure of $U$.
			We start with several properties that follow directly from the definition of $U$-optimality.
			Here the positivity property should not be over-interpreted: it follows from the well-definedness assumption.
			
			\begin{prp}[Properties of optimal e-values]\label{prp:optimizer_properties}
				Let $\varepsilon^*$ be a $U$-optimal e-value.
				\begin{itemize}
					\item \textit{Strict monotonicity}. If $U$ is strictly increasing and $Q_a(\mathcal X) > 0$, then $\textnormal{E}^P[\varepsilon^*] = 1$.
					
					\item \textit{Uniqueness}. If $U$ is strictly concave, then the $U$-optimal e-value is unique up to $Q_a$-almost sure equality.
					
					\item \textit{Positivity}. If $U(0) = -\infty$, then $\varepsilon^* > 0$, $Q_a$-almost surely.
				\end{itemize}
			\end{prp}
			
			In Corollary \ref{cor:structure_pointwise}, we describe properties of the pointwise optimizer \eqref{opt:U-optimal_e-value} in Theorem \ref{thm:EU_NP} that are implied by structure imposed on $U$.
			In the result, $\partial U$ denotes the superdifferential of $U$.
			
			\begin{cor}[Pointwise structure]\label{cor:structure_pointwise}
				Suppose that $\lambda$ and $\varepsilon^*$ satisfy the assumptions of Theorem \ref{thm:EU_NP}.
				\begin{itemize}
					\item \textit{Positive multiplier}. If $U$ is strictly increasing and $Q_a(\mathcal X)>0$, then $\lambda > 0$ and so $\textnormal{E}^P[\varepsilon^*] = 1$.
					
					\item \textit{Supergradient form}. If $U$ is concave, then \eqref{opt:U-optimal_e-value} is equivalent to
					\begin{align*}
						\lambda \Bigg/ \frac{dQ_a}{dP}
							\in \partial U(\varepsilon^*),
						\quad Q_a\textnormal{-a.s.}
					\end{align*}
					
					\item \textit{Positivity}. If $U(0)\in\mathbb R$ and $\lim_{z \downarrow 0} (U(z) - U(0)) / z = \infty$, then $\varepsilon^* > 0$, $Q_a$-almost surely.
					
					\item \textit{Explicit form}. Suppose that $U(0) = \lim_{z \downarrow 0} U(z)$ and that $U$ is differentiable on $(0, \infty)$ with strictly decreasing derivative satisfying $\lim_{z \downarrow 0} U'(z) = \infty$ and $\lim_{z \to \infty} U'(z) = 0$. If $Q_a(\mathcal{X}) > 0$, then $\lambda > 0$ and 
					\begin{align*}
						\varepsilon^*
							= (U')^{-1} \left( \lambda \Bigg/ \frac{dQ_a}{dP} \right), \quad Q_a\textnormal{-a.s.}
					\end{align*}
				\end{itemize}
			\end{cor}

	\section{Optimality for composite hypotheses}
		We now study utility-optimal e-values for a composite hypothesis $H$, which may contain more than one probability measure.
		
		The main insight here is that most of the theory for the simple setting generalizes by replacing $P$ with some utility-dependent measure $P^*$ in an enlargement of $H$.
		Once this measure $P^*$ has been defined, the optimization is the same as in the simple setting, and most of the properties of the optimizer carry over, though some require additional conditions.
				
		\subsection{Decomposition and effective hypotheses}
			To state the results, we require some concepts introduced by \citet{larsson2024numeraire}.
			
			An event is $H$\emph{-negligible} if it has probability zero under every $P \in H$.
			We write
			\begin{align*}
				Q = Q_a + Q_s,
			\end{align*}
			for the generalized Lebesgue decomposition of $Q$ with respect to $H$, where $Q_s$ is concentrated on an $H$-negligible event and $Q_a$ assigns zero mass to every such event.
			If $H$ is simple then this reduces to the classical Lebesgue decomposition.
			
			Using this generalized Lebesgue decomposition, we can extend Definition \ref{dfn:U-optimal} of $U$-optimality to the composite setting without modification.
						
			We denote the collection of valid e-values $\varepsilon : \mathcal{X} \to [0, \infty]$ for $H$ by
			\begin{align*}
				\mathcal{E}(H)
					:= \{\textnormal{measurable }\varepsilon : \textnormal{E}^P[\varepsilon] \leq 1, \textnormal{ for all } P \in H\}.
			\end{align*}
			Notice that every $\varepsilon \in \mathcal{E}(H)$ is finite $Q_a$-almost surely.
			Indeed, $\textnormal{E}^P[\varepsilon] \leq 1$ implies $P(\varepsilon = \infty) = 0$, for every $P \in H$.
			Hence, $\{\varepsilon = \infty\}$ is $H$-negligible and so $Q_a(\varepsilon = \infty) = 0$.
			
			Let $M_+$ denote the collection of nonnegative measures on $\mathcal{X}$.
			The \emph{effective hypothesis} is the largest collection of measures for which every valid e-value for $H$ is valid:
			\begin{align*}
				H^{\textnormal{eff}}
					= \{P \in M_+ : \textnormal{E}^P[\varepsilon] \leq 1, \textnormal{ for all } \varepsilon \in \mathcal{E}(H)\}.
			\end{align*}
			Since the constant e-value $\varepsilon \equiv 1$ is in $\mathcal{E}(H)$, every element of $H^{\textnormal{eff}}$ has total mass at most one.
			This means the effective hypothesis generally contains subprobability measures rather than only probability measures.
			
			\begin{rmk}
				The composite Lebesgue decomposition of $Q$ with respect to $H$ or $H^{\textnormal{eff}}$ is the same, because the hypotheses $H$ and $H^{\textnormal{eff}}$ have the same negligible events.
				Indeed, $H \subseteq H^{\textnormal{eff}}$.
				Conversely, if $N$ is $H$-negligible then $\infty 1_N \in \mathcal{E}(H)$ which forces $P(N) = 0$ for every $P \in H^{\textnormal{eff}}$.
			\end{rmk}

		\subsection{Composite Neyman--Pearson lemma for e-values}
			In Theorem \ref{thm:EU_NP_composite} we present the composite analogue of Theorem \ref{thm:EU_NP}.
			The theorem shows that the composite problem is still a likelihood ratio problem, but with the simple hypothesis $P$ replaced by some element $P^*$ of the effective hypothesis.
			The concavification result in Proposition \ref{prp:concavification} also applies here.
			
			\begin{thm}\label{thm:EU_NP_composite}
				Let $U : [0, \infty] \to [-\infty, \infty]$ be measurable.
				Suppose that there exist $\lambda^* \in [0, \infty)$, $P^* \in H^{\textnormal{eff}}$ with $P^* \ll Q_a$, and a valid e-value $\varepsilon^* \in \mathcal{E}(H)$, satisfying
				\begin{align}\label{opt:composite_U-optimal}
					\varepsilon^*(x)
						\in \argmax_{z \in [0, \infty)} U(z) - \lambda^*\frac{dP^*}{dQ_a}(x) z,\quad Q_a\textnormal{-a.s.},
				\end{align}
				with $\lambda^*(1 - \textnormal{E}^{P^*}[\varepsilon^*]) = 0$, and suppose that $U(\varepsilon^*)$ is $Q_a$-almost surely finite.
				Then $\varepsilon^*$ is $U$-optimal.
			\end{thm}
					
			\begin{rmk}[Recovering Theorem \ref{thm:EU_NP}]
				To recover Theorem \ref{thm:EU_NP} from Theorem \ref{thm:EU_NP_composite}, we set $H = \{P\}$ and choose $P^*$ as the restriction of $P$ to $\{dQ_a/dP > 0\}$.
				Then $P^* \in H^\textnormal{eff}$ and $P^* \ll Q_a$, and $dP^*/dQ_a = 1/(dQ_a/dP)$, $Q_a$-almost surely.
			\end{rmk}

			\begin{rmk}[Why $H^\textnormal{eff}$?]
				Theorem \ref{thm:EU_NP_composite} also holds if $P^*$ is restricted to $H$ or any subset of $H^\textnormal{eff}$.
				However, this reduces the chances that there exists a $P^*$ that satisfies the conditions (see Example \ref{exm:toy_existence}).
		        The effective hypothesis $H^{\textnormal{eff}}$ is the largest class for which the proof goes through, explaining why this object appears.
			\end{rmk}

	        \begin{exm}\label{exm:toy_existence}
	        	Let $\mathcal{X} = \{a, b\}$ and consider $H = \{\delta_a, \delta_b\}$, $Q = .5 \delta_a + .5 \delta_b$ and $U = \log$.
	        	The class of valid e-values is $\mathcal{E}(H) = \{\varepsilon : \varepsilon(a) \leq 1, \varepsilon(b) \leq 1\}$.
	        	Here, the log-optimal e-value is $\varepsilon^* \equiv 1$, but if we were to restrict $P^*$ to $H$ we would not recover this.
	        	Indeed, for $P = \delta_a$ we get $dP/dQ(a) = 2$ and $dP/dQ(b) = 0$.
	        	This means the pointwise optimization \eqref{opt:composite_U-optimal} at $b$ becomes $\argmax_{z \in [0, \infty)} U(z)$.
	        	This has no finite maximizer so that $\varepsilon^*$ is not recovered for $P = \delta_a$, and by symmetry it is not recovered for $P = \delta_b$.
	        	If we instead take $P^* = .5 \delta_a + .5 \delta_b \in H^{\textnormal{eff}}$ then we recover $\varepsilon^*$ with $\lambda^* = 1$.
	        \end{exm}	
			
		\subsection{A converse through a first-order condition}
			As in the simple setting, we also seek a converse to Theorem \ref{thm:EU_NP_composite}.
			There is an important difference between the two settings: for a simple hypothesis, the null distribution $P$ is already known so that the converse only needs to recover the scalar multiplier $\lambda$.
			In the composite setting, the supporting measure $P^*$ is itself unknown.
			We therefore proceed through a first-order condition that allows us to select one measure $P^*$ that supports all feasible directions simultaneously.
			
			For a concave utility function $U$ and valid e-values $\varepsilon^*, \varepsilon \in \mathcal{E}(H)$ such that $U(\varepsilon^*)$ is $Q_a$-almost surely finite, define the pointwise one-sided directional derivative
			\begin{align}\label{dfn:directional_utility}
				\mathrm{D}U(\varepsilon^*; \varepsilon - \varepsilon^*)
					:= \lim_{t \downarrow 0}
					\frac{U((1-t)\varepsilon^* + t\varepsilon) - U(\varepsilon^*)}{t}.
			\end{align}
			The limit exists $Q_a$-almost surely in the extended reals by concavity of $U$.
			
			Now, the first-order condition 
			\begin{align}\label{ineq:integrable_directional_optimality}
				\textnormal{E}^{Q_a}\left[\mathrm{D} U(\varepsilon^*; \varepsilon - \varepsilon^*)\right]
					\leq 0, \textnormal{ for every } \varepsilon \in \mathcal{E}(H),
			\end{align}
			is already sufficient for $U$-optimality, assuming $\mathrm{D}U(\varepsilon^*; \varepsilon-\varepsilon^*)$ is $Q_a$-integrable for every valid e-value $\varepsilon$.
			Indeed, concavity gives, for every valid e-value $\varepsilon$,
			\begin{align*}
				U(\varepsilon) - U(\varepsilon^*)
					\leq \mathrm{D}U(\varepsilon^*; \varepsilon - \varepsilon^*).
			\end{align*}
			
			Proposition \ref{prp:EU_converse_composite} shows that this first-order condition also yields a single supporting measure $P^* \in H^{\textnormal{eff}}$ and the pointwise representation \eqref{opt:composite_U-optimal}.
			
			\begin{prp}[Supporting measure under the first-order condition]\label{prp:EU_converse_composite}
				Assume that $U$ is concave and non-decreasing, and finite on $(0,\infty)$.
				Let $\varepsilon^* \in \mathcal{E}(H)$ be such that $U(\varepsilon^*)$ is $Q_a$-almost surely finite.
				
				Suppose that, for every $\varepsilon \in \mathcal{E}(H)$, $\mathrm{D}U(\varepsilon^*; \varepsilon-\varepsilon^*)$ is $Q_a$-integrable and \eqref{ineq:integrable_directional_optimality} holds.
				Then there exist $\lambda^* \in [0,\infty)$ and $P^* \in H^{\textnormal{eff}}$ with $P^* \ll Q_a$ such that $\varepsilon^*$ satisfies \eqref{opt:composite_U-optimal} and $\lambda^*(1 - \textnormal{E}^{P^*}[\varepsilon^*]) = 0$.
			\end{prp}

			The following corollary now gives the desired converse to Theorem \ref{thm:EU_NP_composite}.
			It starts from a $U$-optimal e-value $\varepsilon^*$ and gives a convenient condition under which the directional first-order condition holds, so that $\varepsilon^*$ can be represented as in \eqref{opt:composite_U-optimal}.
			As a consequence, under condition \eqref{ineq:radial_utility_integrability} the sufficient pointwise condition of Theorem \ref{thm:EU_NP_composite} is also necessary. 
						
			\begin{cor}[Converse to Theorem \ref{thm:EU_NP_composite}]\label{cor:EU_converse_composite_radial}
				Assume that $U$ is concave and non-decreasing, and finite on $(0,\infty)$.
				Suppose that $\varepsilon^* \in \mathcal{E}(H)$ is $U$-optimal and that, for some $a \in (0,1)$,
				\begin{align}\label{ineq:radial_utility_integrability}
					\textnormal{E}^{Q_a}\left[U(\varepsilon^*) - U(a\varepsilon^*)\right]
						< \infty.
				\end{align}
				Then $\varepsilon^*$ satisfies the integrable first-order condition \eqref{ineq:integrable_directional_optimality}.
				As a consequence, there exist $\lambda^* \in [0, \infty)$ and $P^* \in H^\textnormal{eff}$ with $P^* \ll Q_a$ such that $\varepsilon^*$ satisfies \eqref{opt:composite_U-optimal} and $\lambda^*(1 - \textnormal{E}^{P^*}[\varepsilon^*]) = 0$.
			\end{cor}
			
		\subsection{Existence}
			As in the existence results for the simple setting, we distinguish a condition for a fixed pair $(H,Q)$ from a condition on the utility function that guarantees existence for every pair.
			
			For a simple hypothesis, knowledge of $P$ generates a natural ordered family of candidate optimizers $\lambda \mapsto \varepsilon_\lambda$.
			No analogous family is available in advance in the composite setting, because the measure $P^*$ is unknown.
			We therefore impose a uniform integrability condition on a more generic class of candidate optimizers.
			Under a sublinear growth condition, uniform integrability of the candidates suffices.
				
			\begin{thm}[Composite pair-specific existence]\label{thm:existence_composite}
				Suppose that $U$ is concave and upper semicontinuous.
				Assume there exists some benchmark e-value $\varepsilon_0 \in \mathcal{E}(H)$ such that $U(\varepsilon_0)$ is $Q_a$-almost surely finite and 
				\begin{align*}
					\textnormal{E}^{Q_a}[U(\varepsilon) - U(\varepsilon_0)]
				\end{align*}
				is well-defined for every $\varepsilon \in \mathcal{E}(H)$.
				Define the class of candidate e-values
				\begin{align*}
					C(\varepsilon_0) 
						:= \left\{\varepsilon \in \mathcal{E}(H) :  \textnormal{E}^{Q_a}[U(\varepsilon) - U(\varepsilon_0)] \geq 0\right\}.
				\end{align*}
				
				If $\left\{\left[U(\varepsilon) - U(\varepsilon_0)\right]^+ : \varepsilon \in C(\varepsilon_0)\right\}$
				is uniformly integrable under $Q_a$, then a $U$-optimal e-value exists.
				
				If $U$ is finite on $(0, \infty)$, $\limsup_{z \to \infty} U(z)/z \leq 0$ and $U(\varepsilon_0)$ is $Q_a$-integrable, then the uniform integrability follows from $\sup_{\varepsilon \in C(\varepsilon_0)} \textnormal{E}^{Q_a}[\varepsilon] < \infty$.
			\end{thm}
						
			We use the first claim of Theorem \ref{thm:existence_composite} to derive Corollary \ref{cor:log_growth_composite}, which shows that the same log-growth condition on $U$ suffices to guarantee the existence for every pair as in the simple setting.
			\begin{cor}[Composite existence under log-growth]\label{cor:log_growth_composite}
				Under the conditions of Corollary \ref{cor:existence_log_growth}, a $U$-optimal e-value exists.
			\end{cor}
			
			\begin{rmk}
				The proof of Theorem \ref{thm:existence_composite} uses a Koml\'os argument inspired by \citet{larsson2024numeraire}.
				Their existence result follows from Corollary \ref{cor:log_growth_composite}.
				The proof of Corollary \ref{cor:log_growth_composite} follows by constructing a benchmark e-value $\varepsilon_0$ from the log-optimal e-value studied in \citet{larsson2024numeraire}, and subsequently applying Theorem \ref{thm:existence_composite}.
			\end{rmk}

		\subsection{Properties of composite optimizers}
			Effectively all properties of the optimizer transfer to the composite setting.
			The properties from Proposition \ref{prp:optimizer_properties} that depend only on $U$-optimality transfer directly.
			
			\begin{prp}[Properties $U$-optimal composite]\label{prp:composite_optimizer_properties}
				Let $\varepsilon^*$ be a $U$-optimal e-value for $H$.
				Then the Strict monotonicity, Uniqueness and Positivity conclusions of Proposition \ref{prp:optimizer_properties} remain valid, with $\textnormal{E}^P[\varepsilon^*]$ in the Strict monotonicity conclusion replaced by $\sup_{P\in H} \textnormal{E}^P[\varepsilon^*]$.
			\end{prp}

			The pointwise properties from Corollary \ref{cor:structure_pointwise} similarly transfer once a supporting measure $P^*$ has been found.
			In particular, if $U$ is concave then
		    \begin{align*}
		        \lambda^*
		            \frac{dP^*}{dQ_a}
		            \in \partial U(\varepsilon^*),\quad Q_a\textnormal{-a.s.},
		    \end{align*}
		   	while under regularity conditions specified in Corollary \ref{cor:structure_pointwise} we obtain 
		    \begin{align*}
		        \varepsilon^*
		            = (U')^{-1}\left(\lambda^*\frac{dP^*}{dQ_a}\right),\quad Q_a\textnormal{-a.s.}
		    \end{align*} 
			The well-definedness lemmas also generalize.

			\begin{cor}[Pointwise structure composite]
			\label{cor:composite_structure}
				Suppose that $(\lambda^*,P^*,\varepsilon^*)$ satisfies the assumptions of Theorem \ref{thm:EU_NP_composite}.
				Then the conclusions of Corollary \ref{cor:structure_pointwise} remain valid under the substitution of $\lambda / (dQ_a / dP)$ by $\lambda^* (dP^*/dQ_a)$, where $\textnormal{E}^{P^*}[\varepsilon^*] = 1$ and $\lambda^* > 0$ for strictly increasing $U$ if $Q_a(\mathcal{X}) > 0$.
			\end{cor}
			
			\begin{lem}[Composite well-definedness]\label{lem:composite_well_defined}
			    Lemma \ref{lem:well-defined} and Lemma \ref{lem:well-defined_U} remain valid for composite hypotheses.
			\end{lem}

	\section{Applications to Utility functions}
		\subsection{Neyman--Pearson utility}\label{sec:NP_utility}
			We start by considering the utility function that recovers the classical Neyman--Pearson lemma, which we name the `Neyman--Pearson utility function':
			\begin{align*}
				U_\alpha^{\textnormal{NP}}(x)
					= x \wedge 1/\alpha,
			\end{align*}
			for $\alpha \in (0, 1)$, where $x \wedge y = \min(x, y)$.
			
			To understand why this recovers the classical Neyman--Pearson setting, recall that it is classically formulated for a $[0, 1]$-valued test $\tau_\alpha$.
			Following Section \ref{sec:randomized}, $\tau_\alpha$ then maximizes the probability to reject the hypothesis (including randomization) if it maximizes
			\begin{align}\label{obj:NP_01}
				\textnormal{E}^{Q}[\tau_\alpha].
			\end{align}
			On the $H$-negligible event supporting $Q_s$, we can set $\tau_\alpha = 1$ without affecting the validity of the test.
			After doing so, we obtain
			\begin{align*}
				\textnormal{E}^{Q}[\tau_\alpha]
					= Q_s(\mathcal X) + \textnormal{E}^{Q_a}[\tau_\alpha],
			\end{align*}
			so that maximizing classical power corresponds to maximizing $\textnormal{E}^{Q_a}[\tau_\alpha]$.
			
			Rescaling $\varepsilon_\alpha = \tau_\alpha/\alpha$ from $[0, 1]$ to $[0, 1/\alpha]$, this corresponds to maximizing
			\begin{align*}
				\textnormal{E}^{Q_a}[\varepsilon_\alpha]
					= \textnormal{E}^{Q_a}[U_\alpha^{\textnormal{NP}}(\varepsilon_\alpha)].
			\end{align*}
			Finally, as $U_\alpha^{\textnormal{NP}}$ does not reward evidence above $1/\alpha$, we can always choose an optimizer to take values in $[0, 1/\alpha]$.
		
			Combining our results, we obtain the following composite version of the classical Neyman--Pearson lemma, without any assumption on the hypothesis $H$ or alternative $Q$.
			The first-order condition characterization \eqref{ineq:integrable_directional_optimality} collapses into the same statement as $U_\alpha^{\textnormal{NP}}$-optimality here, so we omit it.
			
			\begin{cor}\label{cor:NP}
				A $U_\alpha^{\textnormal{NP}}$-optimal e-value exists for every pair $(H, Q)$ and can be chosen to take values in $[0, 1/\alpha]$.
				For a valid e-value $\varepsilon^{\textnormal{NP}}$ taking values in $[0, 1/\alpha]$, the following statements are equivalent:
				\begin{enumerate}
					\item $\varepsilon^{\textnormal{NP}}$ is $U_\alpha^{\textnormal{NP}}$-optimal,
					\item There exists some $P^* \in H^\textnormal{eff}$ with $P^* \ll Q_a$ and $\lambda^* \in [0, \infty)$ such that, $Q_a$-almost surely,
					\begin{align*}
						\varepsilon^{\textnormal{NP}}(x)
							= 
							\begin{cases}
								1/\alpha, \quad &\textnormal{ if } \lambda^* dP^*/dQ_a(x) < 1, \\
								\kappa_\alpha(x), \quad &\textnormal{ if } \lambda^* dP^*/dQ_a(x) = 1, \\
								0, \quad &\textnormal{ if } \lambda^* dP^*/dQ_a(x) > 1, \\
							\end{cases}
					\end{align*}
					where $\kappa_\alpha : \mathcal{X} \to [0, 1/\alpha]$ is measurable, and we have $\lambda^*(\textnormal{E}^{P^*}[\varepsilon^{\textnormal{NP}}] - 1) = 0$.
				\end{enumerate}
			\end{cor}

			Here, the outcome $\kappa_\alpha$ when $\{\lambda^* dP^*/dQ_a = 1\}$ is classically interpreted as an instruction to reject with external probability $\alpha \kappa_\alpha$.
			We instead interpret $\kappa_\alpha$ directly as some smaller amount of evidence than $1/\alpha$.
			
			\begin{rmk}[Relationship to classical testing]
				An interpretation of the fact that the Neyman--Pearson setting appears for a special utility function means that the classical testing literature has unknowingly been using e-values all along, but under Neyman--Pearson utility.
				We do not believe this utility function matches the true utility of evidence of most practitioners.
				Indeed, we expect that they would usually at least slightly prefer a rejection at, say, level $1\%$ over a rejection at level $5\%$.
				However, Neyman--Pearson utility with $\alpha = 5\%$ states that both are valued the same.
			\end{rmk}
			
			\begin{rmk}[Binary behavior]
				In practice, the event $\lambda^* dP^*/dQ_a = 1$ often has negligible or even zero probability.
				This means the optimal e-value is effectively binary: $\{0, 1/\alpha\}$-valued.
				As a consequence, we must move to a different utility function if we want to move beyond classical binary testing.
				We suspect that implicitly sticking to Neyman--Pearson utility is the reason that $[0, 1]$-valued testing has remained underappreciated, even by \citet{geyer2005fuzzy} who do consider the continuous interpretation of such a test.
			\end{rmk}

			\begin{rmk}[Omitting middle case]\label{rmk:omitting_middle_clause}
				Although the middle case is often omitted from presentations of the standard Neyman--Pearson lemma, it cannot in general be removed from the theorem.
				If the optimization is restricted to $\{0, 1/\alpha\}$-valued tests, a most powerful test need not be a likelihood-ratio threshold test; see Problems 3.17 and 3.18 of \citet{lehmann2022testing}.
			\end{rmk}

		\subsection{Log utility}\label{sec:log_utility}
			The near-universal choice in the literature on e-values is to maximize the expected logarithm under the alternative \citep{grunwald2023safe, larsson2024numeraire}:
			\begin{align*}
				U^{\textnormal{log}}(x) = \log(x).
			\end{align*}
			
			From a utility perspective, it states that the utility increases linearly in the order of magnitude of the evidence.
			For example, a rejection at level $1\%$ receives twice the utility of a rejection at level $10\%$.
			See Remark \ref{rmk:motivating_log_utility} for other motivations.
			
			For log-utility, the $U^{\textnormal{log}}$-optimality of $\varepsilon^*$ becomes
			\begin{align*}
				\textnormal{E}^{Q_a}[\log(\varepsilon / \varepsilon^*)] \leq 0,
			\end{align*}
			for every valid e-value $\varepsilon$, under the convention that $x / 0 = \infty$ for $x > 0$.
			
			Applying our results, we obtain Corollary \ref{cor:log_utility}, which recovers the main findings of \citet{larsson2024numeraire}.
			
			\begin{cor}\label{cor:log_utility}
				Suppose $Q_a(\mathcal{X}) > 0$.
				There exists a $U^\textnormal{log}$-optimal e-value, which is $Q_a$-almost surely unique and strictly positive.
				For a valid e-value $\varepsilon^\textnormal{log}$, the following statements are equivalent:
				\begin{enumerate}
					\item $\varepsilon^\textnormal{log}$ is $U^\textnormal{log}$-optimal,
					\item $\textnormal{E}^{Q_a}[\varepsilon / \varepsilon^\textnormal{log} - 1] \leq 0$, for every $\varepsilon \in \mathcal{E}(H)$,
					\item There exists some $P^* \in H^\textnormal{eff}$ mutually absolutely continuous with $Q_a$ such that
					\begin{align*}
						\varepsilon^\textnormal{log}
							= \frac{1}{Q_a(\mathcal{X})} \frac{dQ_a}{dP^*},\quad Q_a\textnormal{-almost surely,}
					\end{align*} 
					and $\textnormal{E}^{P^*}[\varepsilon^\textnormal{log}] = 1$.
				\end{enumerate}
			\end{cor}
			The factor $1 / Q_a(\mathcal{X})$ appears because $Q_a$ need not be a probability measure.
			Interestingly, the Lagrange multiplier from Theorem \ref{thm:EU_NP_composite} is $\lambda^* = Q_a(\mathcal{X})$ here.

			\begin{rmk}\label{rmk:motivating_log_utility}		
				There are several common motivations for this objective.
				The first motivation is that if $Q_a = Q$ then the optimizer is the likelihood ratio itself, which is well-understood.
				Another common motivation comes from i.i.d. sequential settings, where maximizing this objective at every point in time minimizes the expected time to exceed any threshold (up to overshooting).
				A further motivation is a conditioning-invariance property: the $P$-conditional expectation of a log-utility-optimal e-value given a sub-$\sigma$-algebra is log-utility optimal for the problem restricted to that sub-$\sigma$-algebra.
				In fact, log utility is in some sense the only utility with this property (see Proposition 6 in \citet{koning2025sequentializing}).
			\end{rmk}

		\subsection{Power utility and capping}			
			In this section, we discuss a two-parameter family of utility functions that bridges Neyman--Pearson and log utility.
			
			We start with the family of power utilities, which contains the log utility as a special case
			\begin{align*}
				U_h(z) 
					=
					\begin{cases}
						\frac{z^h - 1}{h}, &\quad h \neq 0, \\
						\log(z), &\quad h = 0,
					\end{cases}
			\end{align*}
			where $\log(z)$ appears as the limit as $h \to 0$.
			For $h \leq 1$, the utility $U_h$ is non-decreasing and concave, and it is strictly concave if $h < 1$.
			
			We now introduce an additional parameter $\alpha \in [0, 1]$, by defining the \emph{capped power utility}
			\begin{align*}
				U_{\alpha, h}(z)
					= U_h(z \wedge 1/\alpha),
			\end{align*}
			where we take $1/0 = \infty$.
			Here, the parameter $h$ controls the curvature of utility, while $\alpha$ determines the largest amount of evidence that receives additional utility.
			This capped power utility class recovers the Neyman--Pearson utility for $h = 1$ and $\alpha \in (0, 1]$ (up to irrelevant constants) and log utility for $h = 0$ and $\alpha = 0$.
			
			In Corollary \ref{cor:capped_power_characterization}, we present the characterization of the $U_{\alpha, h}$-optimal e-values.
			Here, we do not repeat the $\alpha > 0$ and $h = 1$ case already covered by Corollary \ref{cor:NP}.
			The only remaining corner case $(\alpha, h) = (0, 1)$ is discussed in Remark \ref{rmk:uncapped_linear_utility}.
			In the result, we use $1/0=\infty$, $0^r=\infty$ for $r<0$ and $0\times\infty=0$.
			
			\begin{cor}[Capped power optimizers]\label{cor:capped_power_characterization}
				Let $\alpha \in [0, 1]$ and $h < 1$, and let $\varepsilon^{(\alpha,h)}$ be a valid e-value taking values in $[0, 1/\alpha]$.
				If $\alpha = 0$ and $0 < h < 1$, assume additionally that
				\begin{align*}
					\textnormal{E}^{Q_a}\left[\left(\varepsilon^{(\alpha,h)}\right)^h\right]
						< \infty.
				\end{align*}
				
				Then the following statements are equivalent:
				\begin{enumerate}
					\item $\varepsilon^{(\alpha,h)}$ is $U_{\alpha,h}$-optimal,
					\item for every valid e-value $\varepsilon$ taking values in $[0,1/\alpha]$, $\left(\varepsilon^{(\alpha,h)}\right)^{h-1}\left(\varepsilon-\varepsilon^{(\alpha,h)}\right)$ is $Q_a$-integrable and
					\begin{align}\label{ineq:capped_power_first_order}
						\textnormal{E}^{Q_a}\left[\left(\varepsilon^{(\alpha,h)}\right)^{h-1}\left(\varepsilon-\varepsilon^{(\alpha,h)}\right)\right]
							\leq 0,
					\end{align}
					\item there exist a multiplier $\lambda^* \in [0,\infty)$ and a supporting measure $P^*\in H^{\textnormal{eff}}$ with $P^*\ll Q_a$ such that
					\begin{align}\label{opt:capped_power_utility}
						\varepsilon^{(\alpha,h)}
							= \left(\lambda^* \frac{dP^*}{dQ_a}\right)^{-1/(1-h)} \wedge \frac{1}{\alpha}, \quad Q_a\textnormal{-a.s.},
					\end{align}
					and $\lambda^*\left(1-\textnormal{E}^{P^*}[\varepsilon^{(\alpha,h)}]\right) = 0$.
				\end{enumerate}
				Whenever these statements hold, $\varepsilon^{(\alpha,h)} > 0$, $Q_a$-almost surely.
			\end{cor}
		
			We separately present the existence in Corollary \ref{cor:capped_power_existence}, because the existence of a $U$-optimal e-value for every pair $(H, Q)$ is not guaranteed for $\alpha = 0$ and $h > 0$.
			For this parameter range, we rely on our pair-specific existence result Theorem \ref{thm:existence_composite}.
			For $\alpha = 0$ and $0 < h \leq 1$, define
			\begin{align*}
				C_h(\varepsilon_0)
					:= \left\{\varepsilon \in \mathcal{E}(H): \textnormal{E}^{Q_a}[\varepsilon^h - \varepsilon_0^h] \geq 0\right\}.
			\end{align*}
			
			\begin{cor}[Existence for capped power utility]\label{cor:capped_power_existence}
				Let $\alpha \in [0,1]$ and $h \leq 1$.
				We distinguish three cases.
				\begin{enumerate}
					\item If $\alpha > 0$ or $h \leq 0$, then a $U_{\alpha,h}$-optimal e-value exists for every pair $(H, Q)$ and can be chosen to take values in $[0, 1/\alpha]$.
					\item If $\alpha = 0$ and $0 < h < 1$, then a $U_{0, h}$-optimal e-value exists if
						\begin{align*}
							\sup_{\varepsilon \in C_h(\varepsilon_0)} \textnormal{E}^{Q_a}[\varepsilon]
								< \infty,
						\end{align*}
						for some $\varepsilon_0 \in \mathcal{E}(H)$ with $\textnormal{E}^{Q_a}[\varepsilon_0^h] < \infty$.
					\item If $\alpha = 0$ and $h = 1$, then a $U_{0, 1}$-optimal e-value exists if $C_1(\varepsilon_0)$ is uniformly $Q_a$-integrable for some $\varepsilon_0 \in \mathcal{E}(H)$  with $\textnormal{E}^{Q_a}[\varepsilon_0] < \infty$.
				\end{enumerate}
			\end{cor}
			
			\begin{rmk}[$\alpha = 0$ and $h = 1$]\label{rmk:uncapped_linear_utility}
				For $(\alpha, h) = (0, 1)$, we have $U_{0, 1}(z) = z-1$, so the objective is linear and an optimizer need not exist.
				For example, under a simple hypothesis $H = \{P\}$, let $M := \esssup_P \frac{dQ_a}{dP}$.
				If $M < \infty$, an optimizer exists if and only if the likelihood ratio attains $M$ on an event of positive $P$-probability.
				In that case, an optimizer is $\varepsilon^{(0,1)} = 1\{dQ_a/dP=M\} /P\left(\frac{dQ_a}{dP} = M\right)$.
			\end{rmk}
						
			\begin{rmk}
				Uncapped power utility for $h \leq 0$ was also considered by \citet{larsson2024numeraire}, which is exactly the range of $h$ for which existence is guaranteed.
				By restricting to $h \leq 0$ and (implicitly) $\alpha = 0$, this misses the connection to the Neyman--Pearson utility $h = 1$, $\alpha > 0$, and with that our finding that both classical Neyman--Pearson-style testing and log-optimal e-values can be cast into a single expected-utility framework.
			\end{rmk}
			
		\subsection{Step utility}
			As a final example, we consider a step utility which allows a statistician to directly specify their preferences over a collection of conventional evidence thresholds.
			For example, 
			\begin{align*}
				U(x)
					= w_{10} 1\{x \geq 10\} + w_5 1\{x \geq 20\} + w_1 1\{x \geq 100\},
			\end{align*}
			for some weights $w_1, w_5, w_{10} \geq 0$.
			This utility function assigns cumulative preferences to three common rejection levels $1\%$, $5\%$ and $10\%$.
			
			This utility function is notably not concave, such that the converse result Proposition \ref{prp:EU_converse_composite} and existence result Theorem \ref{thm:existence_composite} do not apply.
			Theorem \ref{thm:EU_NP_composite} does apply and provides a way to verify the optimality of a proposed optimizer.
			In particular, given a pair $(P^*, \lambda^*)$, we can write out the finite optimization problem
			\begin{align*}
				\varepsilon^*(x)
					\in \argmax_{z \in \{0, 10, 20, 100\}} U(z) - \lambda^*\frac{dP^*}{dQ_a}(x) z.
			\end{align*}
			If the resulting e-value $\varepsilon^*$ is valid and complementary slackness holds, then Theorem \ref{thm:EU_NP_composite} shows that it is $U$-optimal.
			
			A pragmatic approach to be able to apply our results is to `make $U$ concave', by replacing it with its concave majorant as in Proposition \ref{prp:concavification}.
						
	\section{Discussion}		
		While the optimality theory for simple hypotheses \emph{mathematically} almost entirely carries over to composite hypotheses, the difference has large practical consequences.
		Indeed, replacing the known $P$ by an unknown $P^* \in H^\textnormal{eff}$ that depends on both the alternative $Q$ and utility $U$ makes the computation of an optimal e-value much more challenging.
		
		Another practical consequence is that a single likelihood ratio does not suffice for constructing every utility-optimal e-value in the composite setting: one would need to report an entire collection of likelihood ratios with respect to $P \in H^\textnormal{eff}$.
				
		An open question is to what extent our results extend to arbitrary composite alternatives $\mathcal{Q}$, for a minimax target such as
		\begin{align*}
			\inf_{Q \in \mathcal{Q}} \textnormal{E}^{Q_a}[U(\varepsilon^*) - U(\varepsilon)] \geq 0.
		\end{align*}
		In this context, an interesting observation is that for our relative utility definition of an optimal e-value, the addition of a $Q$-dependent penalty $f(Q)$ cancels, which is sometimes used in non-relative log-utility \citep{grunwald2023safe, RamLarssonRufRamdas2026}
		\begin{align*}
			\inf_{Q \in \mathcal{Q}} &\textnormal{E}^{Q_a}[U(\varepsilon^*) - f(Q) - (U(\varepsilon) - f(Q))] \\
				&= \inf_{Q \in \mathcal{Q}} \textnormal{E}^{Q_a}[U(\varepsilon^*) - U(\varepsilon)].
		\end{align*}
							
	\section{Acknowledgements}
		I thank Alberto Quaini, Olga Kuryatnikova, Ruben van Beesten, Muriel P\'erez-Ortiz, Peter Gr\"unwald, the participants of the CWI e-readers seminar, Nikos Ignatiadis, Aaditya Ramdas, Martin Larsson, Stan Koobs and the participants of the BIRS2025 SAVI workshop.
		
		The author acknowledges the use of ChatGPT-5.6 Pro for proof assistance and proofreading.

		\bibliographystyle{imsart-nameyear} 
		\bibliography{bibliography}       

	\appendix
		
	\section{Proofs}
		\subsection{Proof of Theorem \ref{thm:EU_NP}}
			\begin{proof}
				 Let $\varepsilon$ be some valid e-value.
    			By Lemma~\ref{lem:well-defined}, the relative expected utility is well-defined.
    
				On $\{dQ_a/dP = 0\}$, a comparison to the feasible choice $z = 0$ in the pointwise maximization gives
				\begin{align*}
					- \lambda \varepsilon^* \geq 0,
				\end{align*}
				and hence $\lambda \varepsilon^* = 0$, $P$-almost surely on $\{dQ_a/dP = 0\}$.
				Therefore, by pointwise optimality of $\varepsilon^*$, 
				\begin{align*}
					\textnormal{E}^{Q_a}[U(\varepsilon) - U(\varepsilon^*)]
						&\leq \lambda \textnormal{E}^P[(\varepsilon - \varepsilon^*)1\{dQ_a/dP > 0\}] \\
						&\leq \lambda (\textnormal{E}^P[\varepsilon] - \textnormal{E}^P[\varepsilon^*]) \\
						&= \lambda (\textnormal{E}^P[\varepsilon] - 1)
						\leq 0,
				\end{align*}
				where the equality uses $\lambda (1 - \textnormal{E}^P[\varepsilon^*]) = 0$.
				This proves the result.
			\end{proof}
				
		\subsection{Proof of Lemma \ref{lem:well-defined}}
			\begin{proof}
				We start by assuming all relative utility comparisons are well-defined.
				Taking $\varepsilon = \varepsilon^*$, this means $U(\varepsilon^*) - U(\varepsilon^*)$ must be well-defined $Q_a$-almost surely.
				This requires $U(\varepsilon^*)$ to be finite $Q_a$-almost surely.
				
				For the converse, assume $U(\varepsilon^*)$ is finite $Q_a$-almost surely.
				This implies it is finite $P$-almost surely on $\{dQ_a/dP > 0\}$.
				Pointwise optimality on $\{dQ_a/dP > 0\}$ therefore gives
				\begin{align*}
					\frac{dQ_a}{dP}(U(\varepsilon) - U(\varepsilon^*))
						\leq \lambda (\varepsilon - \varepsilon^*).
				\end{align*}
				As a consequence, its positive part is bounded by $\lambda\varepsilon$:
				\begin{align*}
					\left[\frac{dQ_a}{dP}\left(U(\varepsilon) - U(\varepsilon^*)\right)\right]^+
					\leq \left[\lambda(\varepsilon - \varepsilon^*)\right]^+
						\leq \lambda \varepsilon.
				\end{align*}
				Since $\textnormal{E}^P[\varepsilon] \leq 1$, the positive part is integrable, so that the expected relative utility is well-defined.
			\end{proof}

		\subsection{Proof of Lemma \ref{lem:well-defined_U}}
			\begin{proof}
				As $\varepsilon^*$ is valid, $\varepsilon^* < \infty$ $P$-almost surely, and therefore also $Q_a$-almost surely.
				Since $U(z) < \infty$ for every finite $z$, we therefore also have $U(\varepsilon^*) < \infty$ $Q_a$-almost surely.
				
				Next, let $U(z_0) > -\infty$ for some $z_0 < \infty$.
				Then on $dQ_a / dP > 0$, pointwise optimality of $\varepsilon^*$ gives
				\begin{align*}
					\frac{dQ_a}{dP} U(\varepsilon^*) - \lambda\varepsilon^* 
						\geq \frac{dQ_a}{dP}U(z_0) - \lambda z_0 
						> -\infty.
				\end{align*}
				Since $\varepsilon^* < \infty$ when $dQ_a / dP > 0$, this implies $U(\varepsilon^*) > -\infty$ if $dQ_a / dP > 0$.
				Hence, $U(\varepsilon^*)$ is $Q_a$-almost surely finite.
				
				By Lemma \ref{lem:well-defined}, this is equivalent to well-definedness of the expected relative utility.
			\end{proof}

		\subsection{Proof of Proposition \ref{prp:EU_converse}}
			The result is a pointwise statement about $\varepsilon^*$, whereas $U$-optimality is an aggregate statement, comparing only expected utilities.
			Our proof strategy bridges this gap in two stages.
			
			First, we show that there exists some $\lambda \geq 0$ such that 
			\begin{align}\label{ineq:integrated_lagrange}
				\textnormal{E}^{Q_a}[U(\varepsilon) - U(\varepsilon^*)]
					\leq \lambda (\textnormal{E}^P[\varepsilon] - \textnormal{E}^P[\varepsilon^*]),
			\end{align}
			for every non-negative $\varepsilon$ with finite $P$-expectation for which $U(\varepsilon) - U(\varepsilon^*)$ is $Q_a$-integrable.
			This extends the optimality comparison to competitors that need not be valid e-values.
			
			Second, we apply \eqref{ineq:integrated_lagrange} to competitors that equal a fixed value $z$ on an arbitrary measurable set and equal $\varepsilon^*$ elsewhere.
			This will lead to the desired pointwise maximization claim.

			\begin{proof}
				Since $\varepsilon^*$ is $U$-optimal, all its relative utility comparisons are well-defined by assumption.
				Lemma~\ref{lem:well-defined} therefore gives that $U(\varepsilon^*)$ is finite $Q_a$-almost surely.

				We start with the first stage, in which we construct $\lambda$ according to whether the e-value constraint is binding.
				Throughout this first stage, let $\varepsilon$ be an arbitrary non-negative random variable with finite $P$-expectation for which $U(\varepsilon) - U(\varepsilon^*)$ is $Q_a$-integrable.

				Suppose first that $\textnormal{E}^P[\varepsilon^*] < 1$ and set $\lambda = 0$.
				For sufficiently small $t > 0$, the convex combination
				\begin{align*}
					(1 - t) \varepsilon^* + t\varepsilon,
				\end{align*}
				is valid.
				Concavity and optimality therefore give
				\begin{align*}
					0 
						&\geq \textnormal{E}^{Q_a}[U((1-t) \varepsilon^* + t\varepsilon) - U(\varepsilon^*)] \\
						&\geq t \textnormal{E}^{Q_a}[U(\varepsilon) - U(\varepsilon^*)].
				\end{align*}
				This proves \eqref{ineq:integrated_lagrange} with $\lambda = 0$.
				
				Now suppose that $\textnormal{E}^P[\varepsilon^*] = 1$.
				To find the multiplier $\lambda$, we record the largest relative utility available with a given $P$-budget.
				To this end, we define for $0 \leq t \leq 1$,
				\begin{align*}
					S(t)
						:= \sup_{\eta \geq 0, \textnormal{E}^P[\eta] \leq t} \textnormal{E}^{Q_a}[U(\eta) - U(\varepsilon^*)].
				\end{align*}
				Concavity of $U$ implies that $S$ is concave, while the nesting of the feasible sets implies that it is non-decreasing.
				Moreover, $U$-optimality gives $S(t) \leq 0$ and $S(1) = 0$.
				
				We next verify that $S$ has a finite slope at 1.
				Since $\textnormal{E}^P[\varepsilon^*] = 1$, there exists $n$ such that, for $B := \{1/n \leq \varepsilon^* \leq n\}$,
				\begin{align*}
					\textnormal{E}^P[\varepsilon^*1_B] > 0.
				\end{align*}				
				Define
				\begin{align*}
					\widetilde{\varepsilon}
						:= \varepsilon^*1_{B^c} + \frac{\varepsilon^*}{2} 1_B
						= \varepsilon^* - \frac{\varepsilon^*}{2} 1_B.
				\end{align*}
				Then,
				\begin{align*}
					t_0 
						:= \textnormal{E}^P[\widetilde{\varepsilon}]
						= 1 - \frac{1}{2} \textnormal{E}^P[\varepsilon^*1_B] < 1.
				\end{align*}
				Moreover,
				\begin{align*}
					U(\widetilde{\varepsilon}) - U(\varepsilon^*)
						= (U(\varepsilon^*/2) - U(\varepsilon^*))1_B.
				\end{align*}
				On $B$, we have
				\begin{align*}
					\frac{1}{2n}
						\leq	 \frac{\varepsilon^*}{2}
    					\leq	 \varepsilon^*
    					\leq n.
				\end{align*}
				Now, since $U$ is non-decreasing and finite on $(0,\infty)$,
				\begin{align*}
					U(1 / (2n)) - U(n)
				    \leq U(\varepsilon^* / 2) - U(\varepsilon^*)
				    \leq 0.
				\end{align*}
				This means the relative utility is bounded on $B$, and so has finite $Q_a$-expectation.
				Hence, $S(t_0) > -\infty$, so that $S$ is finite on $[t_0, 1]$.
				
				By concavity, the limit
				\begin{align*}
					\lambda
						:= \lim_{t \uparrow 1} \frac{S(1) - S(t)}{1 - t}
						= \lim_{t \uparrow 1} \frac{-S(t)}{1 - t}
				\end{align*}
				exists.
				It is non-negative, because $S$ is non-decreasing, and finite because, for $t_0 < t < 1$,
				\begin{align*}
					0
						\leq \frac{S(1) - S(t)}{1 - t}
						\leq \frac{S(1) - S(t_0)}{1 - t_0}.
				\end{align*}
				The supporting-line inequality for the concave function $S$ therefore gives
				\begin{align*}
					S(t) 
						\leq \lambda (t - 1),
				\end{align*}
				for all $0 \leq t \leq 1$.
				
				As a consequence, if $\textnormal{E}^P[\varepsilon] \leq 1$, then
				\begin{align*}
					\textnormal{E}^{Q_a}[U(\varepsilon) - U(\varepsilon^*)]
						\leq S(\textnormal{E}^P[\varepsilon])
						\leq \lambda (\textnormal{E}^P[\varepsilon] - 1).
				\end{align*}
				This means \eqref{ineq:integrated_lagrange} holds for all such $\varepsilon$ with $\textnormal{E}^P[\varepsilon] \leq 1$.
				
				We also need \eqref{ineq:integrated_lagrange} for $\varepsilon$ with $\textnormal{E}^P[\varepsilon] > 1$, because the set-wise modifications we use below to prove pointwise optimality need not result in valid e-values.
				Here, the idea is to mix an over-budget $\varepsilon$ with an almost-optimal e-value whose expectation is slightly below one, choosing the mixing weight so that the resulting e-value is valid.
				
				Define $m := \textnormal{E}^P[\varepsilon] > 1$.
				For sufficiently small $s > 0$, choose $\eta_s$ such that $\textnormal{E}^P[\eta_s] \leq 1 - s$ and
				\begin{align*}
					\textnormal{E}^{Q_a}[U(\eta_s) - U(\varepsilon^*)]
						\geq S(1 - s) - s^2.
				\end{align*}
				Set $\theta_s := s / (m - 1 + s)$.
				Then,
				\begin{align*}
					\textnormal{E}^P[(1 - \theta_s) \eta_s + \theta_s \varepsilon]
						\leq (1 - \theta_s) (1 - s) + \theta_s m
						= 1,
				\end{align*}
				so this convex combination of $\eta_s$ and $\varepsilon$ is a valid e-value.
				Concavity and $U$-optimality therefore imply
				\begin{align*}
					0
						&\geq \textnormal{E}^{Q_a}[U((1 - \theta_s)\eta_s + \theta_s \varepsilon) - U(\varepsilon^*)] \\
						&\geq (1 - \theta_s) \textnormal{E}^{Q_a}[U(\eta_s) - U(\varepsilon^*)]
							+ \theta_s \textnormal{E}^{Q_a}[U(\varepsilon) - U(\varepsilon^*)].
				\end{align*}
				Rearranging and using $(1 - \theta_s) / \theta_s = (m - 1) / s$ yields
				\begin{align*}
					\textnormal{E}^{Q_a}[U(\varepsilon) - U(\varepsilon^*)]
						&\leq - \frac{1 - \theta_s}{\theta_s} \textnormal{E}^{Q_a}[U(\eta_s) - U(\varepsilon^*)] \\
						&\leq (m - 1) \left(-\frac{S(1-s)}{s} + s\right).
				\end{align*}
				Now, letting $s \downarrow 0$, we obtain
				\begin{align*}
					\textnormal{E}^{Q_a}[U(\varepsilon) - U(\varepsilon^*)]
						\leq \lambda (m - 1).
				\end{align*}
				This means any utility advantage of an over-budget competitor is bounded by $\lambda$ times its excess $P$-expectation.
				Hence, we also have \eqref{ineq:integrated_lagrange} if $\textnormal{E}^P[\varepsilon] > 1$.
				
				We have now proven \eqref{ineq:integrated_lagrange} for every non-negative $\varepsilon$ with finite $P$-expectation for which $U(\varepsilon) - U(\varepsilon^*)$ is $Q_a$-integrable.
				Moreover, complementary slackness holds: if the constraint is not binding, we choose $\lambda = 0$, while in the binding case $\textnormal{E}^P[\varepsilon^*] = 1$.
				Hence, $\lambda (1 - \textnormal{E}^P[\varepsilon^*]) = 0$.
				
				We now move to the second stage of the proof: applying \eqref{ineq:integrated_lagrange} to competitors $\varepsilon$ that differ from $\varepsilon^*$ only on an arbitrary measurable set.
				
				Fix $z \in (0, \infty)$ and $n \in \mathbb{N}$, and let $A$ be an arbitrary measurable subset of
				\begin{align}\label{set:displayed}
					\left\{\frac{dQ_a}{dP} > 0, \varepsilon^* \leq n, |U(\varepsilon^*)| \leq n\right\}.
				\end{align}
				Define
				\begin{align*}
					\varepsilon
						:= z 1_A + \varepsilon^* 1_{A^c},
				\end{align*}
				which differs from $\varepsilon^*$ only on $A$.
				
				Since $z < \infty$ and $\varepsilon^*$ is valid,
				\begin{align*}
					\textnormal{E}^P[\varepsilon]
						\leq z + \textnormal{E}^P[\varepsilon^*]
						< \infty.
				\end{align*}
				Moreover,
				\begin{align*}
					U(\varepsilon) - U(\varepsilon^*)
						= (U(z) - U(\varepsilon^*)) 1_A,
				\end{align*}
				whose absolute value is bounded by $|U(z)| + n$, and so it is $Q_a$-integrable.
				Applying \eqref{ineq:integrated_lagrange} then gives
				\begin{align*}
					\textnormal{E}^{Q_a}[(U(z) - U(\varepsilon^*)) 1_A]
						\leq \lambda \textnormal{E}^P[(z - \varepsilon^*) 1_A].
				\end{align*}
				Equivalently,
				\begin{align*}
					\textnormal{E}^P\left[1_A \left(\frac{dQ_a}{dP}(U(z) - U(\varepsilon^*)) - \lambda(z - \varepsilon^*)\right)\right]
						\leq 0.
				\end{align*}
				Since $A$ is arbitrary, the expression
				\begin{align*}
					\left(\frac{dQ_a}{dP}(U(z) - U(\varepsilon^*)) - \lambda(z - \varepsilon^*)\right)
				\end{align*}
				must be non-positive $P$-almost surely on the set \eqref{set:displayed}.
				Letting $n \to \infty$ yields
				\begin{align*}
					\frac{dQ_a}{dP} U(\varepsilon^*) - \lambda \varepsilon^*
						\geq \frac{dQ_a}{dP} U(z) - \lambda z,
				\end{align*}
				$P$-almost surely on $\{dQ_a/dP > 0\}$.
				
				Applying the preceding argument to each positive rational $z \in \mathbb Q_{>0}$ and taking the union of the corresponding null sets gives a single $P$-null set outside which
				\begin{align*}
					\frac{dQ_a}{dP}U(\varepsilon^*)-\lambda\varepsilon^*
    					\geq \frac{dQ_a}{dP}U(z)-\lambda z
				\end{align*}
				holds simultaneously for every $z\in\mathbb Q_{>0}$.
				A finite concave function is continuous on the interior of its domain. 
				Hence, by approximating any $z > 0$ by positive rationals, the inequality extends to every $z > 0$.
				Finally, let $z_k \downarrow 0$ with $z_k \in \mathbb Q_{>0}$. 
				Since $U$ is non-decreasing,
				\begin{align*}
					 U(z_k)\geq U(0),
				\end{align*}
				and therefore
				\begin{align*}
					\frac{dQ_a}{dP}U(\varepsilon^*)-\lambda\varepsilon^*
				    \geq \frac{dQ_a}{dP}U(0)-\lambda z_k.
				\end{align*}
				Letting $k \to \infty$ proves the inequality also for $z = 0$.
				
				It remains to verify the pointwise condition on the event $\{dQ_a/dP = 0\}$.
				Consider
				\begin{align*}
					\varepsilon
						:= \varepsilon^* 1\left\{\frac{dQ_a}{dP} > 0\right\}.
				\end{align*}
				This modification does not change $\varepsilon^*$, $Q_a$-almost surely, so the left-hand side of \eqref{ineq:integrated_lagrange} is zero.
				Hence,
				\begin{align*}
					0
						\leq - \lambda \textnormal{E}^P\left[\varepsilon^* 1\left\{\frac{dQ_a}{dP} = 0\right\}\right].
				\end{align*}
				The right-hand side is non-positive, and therefore
				\begin{align*}
					\lambda \varepsilon^* = 0,
				\end{align*}
				$P$-almost surely on $\{dQ_a/dP = 0\}$.
				On this set, the pointwise objective is simply $-\lambda z$.
				Its maximum is zero, and attained by $\varepsilon^*$ because $\lambda\varepsilon^* = 0$.
				
				We have therefore shown that
				\begin{align*}
					\varepsilon^*(x)
						\in \argmax_{z \in [0, \infty)} \frac{dQ_a}{dP}(x) U(z) - \lambda z,
				\end{align*}
				$P$-almost surely, completing the proof.
			\end{proof}

		\subsection{Proof of Proposition \ref{prp:concavification}}
			\begin{proof}
				If $a = 0$, the two objective functions both equal $-bz$, so both conclusions follow.
				Suppose therefore that $a > 0$ and define
				\begin{align*}
					M := \sup_{z \in [0,\infty)} \left\{aU(z) - bz\right\},
				\end{align*}
				and
				\begin{align*}
					\overline{M} := \sup_{z \in [0,\infty)} \left\{a\overline{U}(z) - bz\right\}.
				\end{align*}
				Since $\overline{U} \geq U$, we immediately have
				\begin{align}\label{ineq:concavification_lower_bound}
					\overline{M} \geq M.
				\end{align}
				
				If $M = \infty$, then \eqref{ineq:concavification_lower_bound} gives $\overline{M} = \infty$.
				If $M = -\infty$, then $U(z) = -\infty$ for every $z \in [0, \infty)$, and therefore $\overline{U}(z) = -\infty$ for every $z \in [0,\infty)$, so that $\overline{M} = M$.
				
				It remains to consider the case where $M$ is finite.
				By the definition of $M$,
				\begin{align*}
					aU(z)-bz \leq M,
				\end{align*}
				for every $z \in [0, \infty)$.
				Equivalently,
				\begin{align*}
					U(z) \leq \frac{M+bz}{a}.
				\end{align*}
				Hence, the affine function
				\begin{align*}
					z \mapsto \frac{M+bz}{a}
				\end{align*}
				is a concave majorant of $U$.
				Since $\overline{U}$ is the least concave majorant of $U$,
				\begin{align*}
					\overline{U}(z) \leq \frac{M+bz}{a}.
				\end{align*}
				It follows that
				\begin{align*}
					a\overline{U}(z) - bz \leq M,
				\end{align*}
				for every $z \in [0,\infty)$, and hence
				\begin{align*}
					\overline{M} \leq M.
				\end{align*}
				Together with \eqref{ineq:concavification_lower_bound}, this proves that $\overline{M}=M$.
				
				Finally, let $z^*$ be a maximizer of the original problem.
				Then
				\begin{align*}
					M = aU(z^*)-bz^*
						\leq a\overline{U}(z^*)-bz^*
						\leq \overline{M}
						= M.
				\end{align*}
				All inequalities are therefore equalities, so $z^*$ also maximizes the concavified problem.
			\end{proof}
			
		\subsection{Proof of Theorem \ref{thm:existence_simple}}
			The result separates the two issues that arise in proving existence.
    		We first show that the pointwise maximizers exist and can be selected measurably.
    		Assuming that these pointwise maximizers are $P$-integrable, we then vary $\lambda$ until their $P$-expectation reaches the e-value validity constraint.
    		
    		\begin{proof}
    			Fix $\lambda > 0$.
    			We first verify that pointwise maximizers exist.
    			
    			For a fixed value $\ell \geq 0$, consider the pointwise objective
    			\begin{align*}
    				z \mapsto \ell U(z) - \lambda z,
    			\end{align*}	
    			$z \geq 0$.
    			If $\ell = 0$, this objective equals $-\lambda z$, so its unique maximizer is zero.
    			If $\ell > 0$, then \eqref{ineq:sublinear} gives
    			\begin{align*}
    				\limsup_{z \to \infty} \frac{\ell U(z) - \lambda z}{z}
    					\leq -\lambda
    					< 0.
    			\end{align*}
    			Hence,
    			\begin{align*}
    				\ell U(z) - \lambda z
    					\to -\infty,
    			\end{align*}
    			as $z \to \infty$.
    			Upper semicontinuity therefore gives attainment, while concavity implies that the set of maximizers is a nonempty compact interval.
    			
    			We next verify that the pointwise maximizers can be selected measurably.
    			For every $\ell \geq 0$, choose the smallest maximizer of
    			\begin{align*}
    				z \mapsto \ell U(z) - \lambda z.
    			\end{align*}
    			It suffices to show that this choice is non-decreasing in $\ell$, since every monotone real-valued function is Borel measurable.
    			
    			To this end, let $0 < \ell_1 < \ell_2$, and let $z_i$ be the smallest maximizer corresponding to $\ell_i$.
    			Equivalently, $z_i$ maximizes
    			\begin{align*}
    				z \mapsto U(z) - \frac{\lambda}{\ell_i} z.
    			\end{align*}
    			Pointwise optimality therefore gives
    			\begin{align*}
    				U(z_1) - \frac{\lambda}{\ell_1} z_1
    					&\geq U(z_2) - \frac{\lambda}{\ell_1} z_2, \\
    				U(z_2) - \frac{\lambda}{\ell_2} z_2
    					&\geq U(z_1) - \frac{\lambda}{\ell_2} z_1.
    			\end{align*}
    			Adding these inequalities yields
    			\begin{align*}
    				\lambda \left(\frac{1}{\ell_1} - \frac{1}{\ell_2}\right) (z_2 - z_1)
    					\geq 0,
    			\end{align*}
    			and so $z_1 \leq z_2$.
    			For $\ell = 0$, the unique maximizer is zero.
    			Combining these observations, the smallest maximizer is a non-decreasing, and hence Borel measurable, function of $\ell$.
    			Evaluating this function at $dQ_a / dP$ defines the measurable random variable $\varepsilon_\lambda$.
    			
    			Next, we show that this maximizer $\varepsilon_\lambda$ is non-increasing in $\lambda$.
    			Indeed, let $0 < \lambda_1 < \lambda_2$, and let $z_i$ be the selected pointwise maximizer corresponding to $\lambda_i$.
    			Pointwise optimality then gives
    			\begin{align*}
    				\frac{dQ_a}{dP} U(z_1) - \lambda_1 z_1
    					&\geq \frac{dQ_a}{dP} U(z_2) - \lambda_1 z_2, \\
    				\frac{dQ_a}{dP} U(z_2) - \lambda_2 z_2
    					&\geq \frac{dQ_a}{dP} U(z_1) - \lambda_2 z_1.
    			\end{align*}
    			Adding these inequalities yields
    			\begin{align*}
    				(\lambda_2 - \lambda_1) (z_1 - z_2)
    					\geq 0,
    			\end{align*}
    			and hence $z_1 \geq z_2$.
    			
    			We now assume that
    			\begin{align*}
    				\textnormal{E}^P[\varepsilon_\lambda]
    					< \infty,
    			\end{align*}
    			for every $\lambda > 0$, and show that $\lambda$ can be chosen so that the pointwise maximizer satisfies the e-value validity constraint.
    			
    			We first show that
    			\begin{align}\label{limit:e_lambda_zero}
    				\varepsilon_\lambda \downarrow 0,
    			\end{align}
    			as $\lambda \to \infty$.
    			Monotonicity ensures that the pointwise limit exists.
    			Fix $x$ such that $dQ_a / dP(x) < \infty$ and suppose, for the sake of contradiction, that
    			\begin{align*}
    				\varepsilon_\lambda(x) \downarrow c,
    			\end{align*}
    			for some $c > 0$.
    			We then necessarily have that $dQ_a / dP(x) > 0$, since $\varepsilon_\lambda(x) = 0$ if $dQ_a/dP(x) = 0$.
    			Fix some $\lambda_0 > 0$.
    			For every $\lambda \geq \lambda_0$, monotonicity and convergence to $c$ then give
    			\begin{align*}
    				c 
    					\leq \varepsilon_\lambda(x)
    					\leq \varepsilon_{\lambda_0}(x).
    			\end{align*}
    			Comparison in the pointwise maximization with $c / 2$ therefore gives
    			\begin{align*}
    				\lambda \left(\varepsilon_\lambda(x) - \frac{c}{2}\right)
    					\leq \frac{dQ_a}{dP} \left(U(\varepsilon_\lambda(x)) - U(c/2)\right).
    			\end{align*}
    			Both arguments of $U$ on the right-hand side belong to the fixed compact interval $[\frac{c}{2}, \varepsilon_{\lambda_0}(x)] \subset (0, \infty)$.
    			Since $U$ is finite and concave on $(0, \infty)$, it is bounded on this interval.
    			Hence, the right-hand side remains bounded as $\lambda \to \infty$.
    			On the other hand,
    			\begin{align*}
    				\lambda \left(\varepsilon_\lambda(x) - \frac{c}{2}\right)
    					\geq \frac{\lambda c}{2} \to \infty,
    			\end{align*}
    			which is a contradiction.
    			Hence, the pointwise limit must be zero, proving \eqref{limit:e_lambda_zero}.
    			
    			We have shown that increasing $\lambda$ pushes the pointwise maximizer towards zero. 
    			We now translate this pointwise convergence into convergence of its $P$-expectation.
    			This will ensure that the e-value constraint is satisfied for all sufficiently large $\lambda$.
    			
    			Fix $\lambda_0 > 0$.
    			Since $\varepsilon_\lambda$ is non-increasing in $\lambda$,
 				\begin{align*}
 					0
 						\leq \varepsilon_\lambda
 						\leq \varepsilon_{\lambda_0},
 				\end{align*}
 				for $\lambda \geq \lambda_0$.
 				The random variable $\varepsilon_{\lambda_0}$ is $P$-integrable by assumption.
 				Hence, dominated convergence gives
 				\begin{align}\label{limit:expected_e_lambda_zero}
 					\textnormal{E}^P[\varepsilon_\lambda]
 						\to 0,
 				\end{align}
 				as $\lambda \to \infty$.
 				In particular, $\varepsilon_\lambda$ is a valid e-value for all sufficiently large $\lambda$.

 				It remains to select one multiplier $\lambda$ and corresponding pointwise maximizer $\varepsilon_\lambda$ that satisfy both the e-value constraint and the complementary slackness.
 				Define
 				\begin{align*}
 					\lambda^*
 						:= \inf\left\{\lambda > 0 : \textnormal{E}^P[\varepsilon_\lambda] \leq 1\right\}.
 				\end{align*}
 				This set is non-empty by \eqref{limit:expected_e_lambda_zero}.
 				
 				Since $\varepsilon_\lambda$ is non-increasing in $\lambda$,
 				\begin{align}\label{ineq:threshold_expectations}
 					\textnormal{E}^P[\varepsilon_\lambda]
 						&\leq 1, \quad \textnormal{ for } \lambda > \lambda^*, \\
 					\textnormal{E}^P[\varepsilon_\lambda]
 						&> 1, \quad \textnormal{ for } 0 < \lambda < \lambda^*.
 				\end{align}
				
				We now distinguish whether $\lambda^*$ is positive or zero.
				If $\lambda^* > 0$, the complementary slackness condition means we must construct a pointwise maximizer whose $P$-expectation is exactly 1.
				If $\lambda^* = 0$, the complementary slackness condition holds automatically, so it suffices to find some valid pointwise maximizer.
				
				Suppose first that $\lambda^* > 0$.
				We cannot simply take the selected pointwise maximizer at $\lambda^*$: 
				\begin{align*}
    				\lambda\mapsto\mathbb E^P[\varepsilon_\lambda]
				\end{align*}
				may jump from above one to below one at $\lambda^*$.
				We therefore consider the pointwise limits from above and below:
				\begin{align*}
					\varepsilon^-
						&:= \lim_{\lambda \downarrow \lambda^*, \lambda > \lambda^*} \varepsilon_\lambda, \\
					\varepsilon^+
						&:= \lim_{\lambda \uparrow \lambda^*, \lambda < \lambda^*} \varepsilon_\lambda.
				\end{align*}
				These limits exist by monotonicity, and
				\begin{align*}
					\varepsilon^- \leq \varepsilon^+.
				\end{align*}
				
				We first verify that both limits are indeed pointwise maximizers at $\lambda^*$.
				On $\{dQ_a/dP = 0\}$, all $\varepsilon_\lambda$ equal zero, so this is immediate.
				Now, fix $x$ such that $dQ_a/dP(x) > 0$, and let $\lambda_k \to \lambda^*$ from either side and write
				\begin{align*}
					\varepsilon_{\lambda_k}(x)
						\to \overline{\varepsilon}(x).
				\end{align*}
				For every $z \geq 0$, pointwise optimality gives
				\begin{align*}
					\frac{dQ_a}{dP}(x) U(\varepsilon_{\lambda_k}(x)) - \lambda_k \varepsilon_{\lambda_k}(x)
						\geq \frac{dQ_a}{dP}(x) U(z) - \lambda_k z.
				\end{align*}
				By upper semicontinuity of $U$,
				\begin{align*}
					&\frac{dQ_a}{dP}(x) U(\overline{\varepsilon}(x)) - \lambda^*\overline{\varepsilon}(x) \\
						&\geq \limsup_{k\to\infty} \left(\frac{dQ_a}{dP}(x) U(\varepsilon_{\lambda_k}(x)) - \lambda_k \varepsilon_{\lambda_k}(x)\right) \\
						&\geq \frac{dQ_a}{dP}(x) U(z) - \lambda^* z.
				\end{align*}
				As a consequence, both $\varepsilon^-$ and $\varepsilon^+$ are pointwise maximizers at $\lambda^*$.
 				
 				Next, we show that their expectations lie on opposite sides of 1.
 				As $\lambda \downarrow \lambda^*$ from above, $\varepsilon_\lambda \uparrow \varepsilon^-$, while \eqref{ineq:threshold_expectations} gives $\textnormal{E}^P[\varepsilon_\lambda] \leq 1$.
 				Monotone convergence therefore gives
 				\begin{align*}
 					\textnormal{E}^P[\varepsilon^-] \leq 1.
 				\end{align*}
 				As $\lambda \uparrow \lambda^*$ from below, $\varepsilon_\lambda \downarrow \varepsilon^+$.
 				For $\lambda$ sufficiently close to $\lambda^*$,
 				\begin{align*}
 					0 \leq \varepsilon_\lambda \leq \varepsilon_{\lambda^*/2}.
 				\end{align*}
 				Since $\varepsilon_{\lambda^*/2}$ is $P$-integrable by assumption, both one-sided limits are finite $P$-almost surely.
 				Hence, dominated convergence and \eqref{ineq:threshold_expectations} give
 				\begin{align*}
 					\textnormal{E}^P[\varepsilon^+] \geq 1.
 				\end{align*}
 				As a consequence,
 				\begin{align*}
 					\textnormal{E}^P[\varepsilon^-] \leq 1 \leq \textnormal{E}^P[\varepsilon^+].
 				\end{align*}
 				
 				By linearity of expectations, we may choose $t \in [0, 1]$ such that
 				\begin{align*}
 					\varepsilon^*
 						:= (1-t) \varepsilon^- + t \varepsilon^+
 				\end{align*}
 				satisfies $\textnormal{E}^P[\varepsilon^*] = 1$.
 				For every $x$, the pointwise objective is concave in $z$, so its set of maximizers is an interval.
 				Since both $\varepsilon^-(x)$ and $\varepsilon^+(x)$ are pointwise maximizers at $\lambda^*$, their convex combination $\varepsilon^*(x)$ is also a pointwise maximizer at $\lambda^*$.
				Hence, $\varepsilon^*$ is valid and satisfies
 				\begin{align*}
 					\lambda^* \left(1 - \textnormal{E}^P[\varepsilon^*]\right)
 						= 0.
 				\end{align*}
 				
 				Suppose now that $\lambda^* = 0$.
 				In this case, \eqref{ineq:threshold_expectations} gives
 				\begin{align*}
 					\textnormal{E}^P[\varepsilon_\lambda] \leq 1,
 				\end{align*}
 				for $\lambda > 0$.
 				Hence, every $\varepsilon_\lambda$ is already a valid e-value for $\lambda > 0$.
 				Since $\varepsilon_\lambda$ increases as $\lambda \downarrow 0$, define
 				\begin{align*}
 					\varepsilon^*
 						:= \lim_{\lambda \downarrow 0} \varepsilon_\lambda.
 				\end{align*}
 				Monotone convergence gives
 				\begin{align*}
 					\textnormal{E}^P[\varepsilon^*]
 						= \lim_{\lambda \downarrow 0} \textnormal{E}^P[\varepsilon_\lambda]
 						\leq 1.
 				\end{align*}
 				Hence, $\varepsilon^*$ is valid and so finite $P$-almost surely.
 				
 				We now verify that $\varepsilon^*$ is pointwise optimal at the multiplier zero.
 				On $\{dQ_a/dP > 0\}$, pointwise optimality of $\varepsilon_\lambda$ gives, for every $z \geq 0$,
 				\begin{align*}
 					U(\varepsilon_\lambda)
 						\geq U(z) + \frac{\lambda}{dQ_a/dP} (\varepsilon_\lambda - z).
 				\end{align*}
 				Since $\varepsilon_\lambda \uparrow \varepsilon^* < \infty$, we have
 				\begin{align*}
 					-z \leq \varepsilon_\lambda - z \leq \varepsilon^*.
 				\end{align*}
 				Hence
 				\begin{align*}
 					\frac{\lambda}{dQ_a / dP}(\varepsilon_\lambda - z)
 						\to 0,
 				\end{align*}
 				as $\lambda \downarrow 0$.
 				As a consequence,
 				\begin{align*}
 					\liminf_{\lambda \downarrow 0} U(\varepsilon_\lambda)
 						\geq U(z).
 				\end{align*}
 				Upper semicontinuity of $U$ and $\varepsilon_\lambda \to \varepsilon^*$ then give
 				\begin{align*}
 					U(\varepsilon^*)
 						\geq \limsup_{\lambda \downarrow 0} U(\varepsilon_\lambda)
 						\geq \liminf_{\lambda\downarrow0}U(\varepsilon_\lambda) 
 						\geq U(z).
 				\end{align*}
 				As a consequence, $\varepsilon^*$ is a pointwise maximizer on $\{dQ_a/dP > 0\}$.
 				On $\{dQ_a / dP = 0\}$, the pointwise objective at $\lambda^* = 0$ equals zero, so every $z \geq 0$ is a maximizer.
 				Hence, \eqref{opt:U-optimal_e-value} holds, while complementary slackness is automatic because $\lambda^* = 0$.
 				
 				In both the $\lambda^* > 0$ and $\lambda^* = 0$ cases, we have now constructed a valid pointwise maximizer that satisfies complementary slackness.
 				It remains to verify the finiteness condition in Theorem \ref{thm:EU_NP}.
 				
 				First, concavity and finiteness of $U$ on $(0, \infty)$ rule out $U(0) = \infty$.
 				Indeed,
 				\begin{align*}
 					U(1)
 						\geq \frac{1}{2} U(0) + \frac{1}{2} U(2),
 				\end{align*}
 				which would contradict the finiteness of $U(1)$ if $U(0) = \infty$.
 				Since $\varepsilon^* < \infty$ $P$-almost surely, we therefore have
 				\begin{align*}
 					U(\varepsilon^*) < \infty, \quad P\textnormal{-a.s.}
 				\end{align*}
 				
 				On $\{dQ_a/dP > 0\}$, comparison in the pointwise maximization with $z = 1$ gives
 				\begin{align*}
 					\frac{dQ_a}{dP} U(\varepsilon^*) - \lambda^*\varepsilon^*
 						\geq \frac{dQ_a}{dP} U(1) - \lambda^*
 						> -\infty.
 				\end{align*}
 				It follows that $U(\varepsilon^*) > -\infty$ on $\{dQ_a/dP > 0\}$.
 				Since $Q_a$ is concentrated on this event, $U(\varepsilon^*)$ is finite $Q_a$-almost surely.
 				
 				The pair $(\lambda^*, \varepsilon^*)$ therefore satisfies all the conditions of Theorem \ref{thm:EU_NP}, completing the proof.
    		\end{proof}			
			
		\subsection{Proof of Corollary \ref{cor:existence_log_growth}}
			\begin{proof}
				We need to verify the conditions of Theorem \ref{thm:existence_simple}.
				Taking $x = 1$ in \eqref{ineq:log_growth} yields
				\begin{align*}
					U(z)
						\leq U(1) + K \log(z),
				\end{align*}
				for $z \geq 1$.
				Hence,
				\begin{align*}
					\limsup_{z \to \infty} \frac{U(z)}{z} \leq 0,
				\end{align*}
				so the measurable pointwise maximizers $\varepsilon_\lambda$ as in Theorem \ref{thm:existence_simple} exist.
				
				It remains to verify their $P$-integrability.
				If $\varepsilon_\lambda < 2$, it is bounded by 2.
				If $\varepsilon_\lambda \geq 2$, then a comparison in the pointwise maximization with $\varepsilon_\lambda / 2$ gives
				\begin{align*}
					\frac{\lambda}{2} \varepsilon_\lambda
						&\leq \frac{dQ_a}{dP} \left(U(\varepsilon_\lambda) - U(\varepsilon_\lambda / 2)\right) \\
						&\leq K \log(2) \frac{dQ_a}{dP}.
				\end{align*}
				Hence, in either case,
				\begin{align*}
					\varepsilon_\lambda
						\leq 2 + \frac{2K \log(2)}{\lambda} \frac{dQ_a}{dP}.
				\end{align*}
				Since $\textnormal{E}^P[dQ_a/dP] = Q_a(\mathcal{X}) \leq 1$, it follows that $\textnormal{E}^P[\varepsilon_\lambda] < \infty$ for every $\lambda > 0$.
				Theorem \ref{thm:existence_simple} then proves the result.
			\end{proof}
		
		\subsection{Proof of Proposition \ref{prp:optimizer_properties}}
			\begin{proof}
				We prove the conclusions in order. \\
				
				\noindent\textit{Strict monotonicity}.
				Suppose $\textnormal{E}^P[\varepsilon^*] < 1$.
				Then
				\begin{align*}
					\widetilde{\varepsilon}
						:= \varepsilon^*+\frac{1 - \textnormal{E}^P[\varepsilon^*] }{2}
				\end{align*}
				is a valid e-value.
				Strict monotonicity gives
				\begin{align*}
					U(\widetilde{\varepsilon})
						> U(\varepsilon^*),
				\end{align*}
				$Q_a$-almost surely, contradicting $U$-optimality because $Q_a(\mathcal X)>0$. \\
				
				\noindent\textit{Uniqueness}. 
				Let $\widetilde{\varepsilon}$ be another $U$-optimal e-value.
				Optimality in both directions implies
				\begin{align*}
					\textnormal{E}^{Q_a}\left[U(\widetilde{\varepsilon}) - U(\varepsilon^*)\right]
						=0.
				\end{align*}
				The average
				\begin{align*}
					\overline{\varepsilon} 
						:= \frac{\varepsilon^* + \widetilde{\varepsilon}}{2}
				\end{align*}
				is valid.
				Strict concavity gives
				\begin{align*}
					U(\overline{\varepsilon})
						\geq \frac{1}{2}U(\varepsilon^*) + \frac{1}{2}U(\widetilde{\varepsilon}),
				\end{align*}
				with strict inequality on $\{\widetilde{\varepsilon} \neq \varepsilon^*\}$.
				Rearranging yields
				\begin{align*}
					U(\overline{\varepsilon})-U(\varepsilon^*)
						\geq \frac{1}{2} \left[U(\widetilde{\varepsilon}) - U(\varepsilon^*)\right],
				\end{align*}
				which is strict on $\{\widetilde{\varepsilon} \neq \varepsilon^*\}$.
				Since the right-hand side has expectation zero, we would obtain
				\begin{align*}
					\textnormal{E}^{Q_a}\left[U(\overline{\varepsilon}) - U(\varepsilon^*)\right]
						> 0,
				\end{align*}
				if $Q_a(\widetilde{\varepsilon} \neq \varepsilon^*) > 0$.
				This contradicts $U$-optimality of $\varepsilon^*$, and so $\widetilde{\varepsilon} = \varepsilon^*$, $Q_a$-almost surely \\
				
				\noindent\textit{Positivity}.
				Lemma \ref{lem:well-defined} gives that $U(\varepsilon^*)$ is finite $Q_a$-almost surely.
				Since $U(0)=-\infty$, this requires $\varepsilon^*>0$, $Q_a$-almost surely.
			\end{proof}

		\subsection{Proof of Corollary \ref{cor:structure_pointwise}}
			\begin{proof}
				We again prove the conclusions in order. \\
				
				\noindent\textit{Positive multiplier}.
				If $\lambda = 0$, then on $\{dQ_a/dP > 0\}$ the pointwise condition requires the finite e-value $\varepsilon^*$ to maximize the strictly increasing function $U$ over $[0,\infty)$.
				This is impossible when $Q_a(\mathcal X) > 0$. \\
				
				\noindent\textit{Supergradient form}.
				On $\{dQ_a/dP>0\}$, the pointwise condition is equivalent to
				\begin{align*}
					U(z)-U(\varepsilon^*)
						\leq
					\frac{\lambda}{dQ_a/dP} (z - \varepsilon^*),
				\end{align*}
				for $z \geq 0$.
				This is precisely the definition of the supergradient of a concave function. \\
				
				\noindent\textit{Positivity}.
				Suppose that $\varepsilon^* = 0$ and $dQ_a/dP > 0$. 
				Pointwise optimality then gives, for every $z > 0$,
				\begin{align*}
					\frac{U(z) - U(0)}{z}
						\leq \frac{\lambda}{dQ_a/dP}.
				\end{align*}
				Letting $z \downarrow 0$ gives a contradiction.
				$Q_a$-almost sure positivity follows because $Q_a(dQ_a/dP > 0) = 1$. \\
				
				\noindent\textit{Explicit form}.
				The assumptions imply that $U$ is strictly concave and strictly increasing on $(0, \infty)$, while the right-continuity condition extends concavity to zero.
				The Positive multiplier conclusion therefore gives $\lambda > 0$.
				
				The preceding positivity results also give $\varepsilon^* > 0$, $Q_a$-almost surely.
				Indeed, this follows from Proposition \ref{prp:optimizer_properties} if $U(0) = -\infty$.
				If $U(0)$ is finite, concavity gives
				\begin{align*}
					\frac{U(z) - U(0)}{z}
						\geq U'(z),
				\end{align*}
				so the positivity conclusion above applies.
				
				The Supergradient form conclusion and differentiability now give
				\begin{align*}
					U'(\varepsilon^*)
						=
					\lambda
					\Bigg/
					\frac{dQ_a}{dP}.
				\end{align*}
				Since $U'$ is strictly decreasing with range $(0,\infty)$, applying its inverse proves the result.
			\end{proof}

		\subsection{Proof of Theorem \ref{thm:EU_NP_composite}}	
			\begin{proof}
				Let $\varepsilon \in \mathcal{E}(H)$ be some valid e-value.
   	    		By Lemma~\ref{lem:composite_well_defined}, the relative expected utility is well-defined.
   	    		
   	    		Pointwise optimality gives
   	    		\begin{align*}
   	    			U(\varepsilon) - U(\varepsilon^*)
   	    				\leq \lambda^* \frac{dP^*}{dQ_a} (\varepsilon - \varepsilon^*),
   	    		\end{align*}
   	    		$Q_a$-almost surely.
   	    		Integrating with respect to $Q_a$ yields
   	    		\begin{align*}
   	    			\textnormal{E}^{Q_a}[U(\varepsilon) - U(\varepsilon^*)]
   	    				&\leq \lambda^* \left(\textnormal{E}^{P^*}[\varepsilon] - \textnormal{E}^{P^*}[\varepsilon^*]\right) \\
   	    				&\leq \lambda^* \left(1 - \textnormal{E}^{P^*}[\varepsilon^*]\right)
   	    				= 0,
   	    		\end{align*}
   	    		where the inequality uses $P^* \in H^{\textnormal{eff}}$ and the equality uses complementary slackness.
   	    		Hence, $\varepsilon^*$ is $U$-optimal.
			\end{proof}

		\subsection{Proof of Proposition \ref{prp:EU_converse_composite}}
			We write $L^1(Q_a)$ for the $Q_a$-integrable random variables and $L^\infty(Q_a)$ for the $Q_a$-essentially bounded random variables.
			
			Before we present the proof, we introduce some concepts relating to derivatives of $U$, where we carefully handle what happens near the boundary $0$.
			For every finite $x \geq 0$, define the right derivative
			\begin{align*}
				U'_+(x)
					:= \lim_{h \downarrow 0}
					\frac{U(x + h) - U(x)}{h},
			\end{align*}
			and, for every finite $x > 0$, define the left derivative
			\begin{align*}
				U'_-(x)
					:= \lim_{h \downarrow 0}
					\frac{U(x) - U(x - h)}{h}.
			\end{align*}
			The right derivative at zero is understood in the extended reals.
			In particular, if $U(0) = -\infty$, then $U'_+(0) = \infty$ since $U$ is finite on $(0, \infty)$.
			The case $U(0) = \infty$ is ruled out by concavity and finiteness of $U$ on $(0,\infty)$.
			
			For $x>0$, concavity, monotonicity and finiteness of $U$ imply
			\begin{align*}
				0 	\leq U'_+(x)
					\leq U'_-(x)
					< \infty,
			\end{align*}
			and
			\begin{align*}
				\partial U(x)
					= [U'_+(x),U'_-(x)].
			\end{align*}
			Moreover, the functions $U'_+$ and $U'_-$ are non-increasing on $(0,\infty)$, and hence Borel measurable.
			
			We introduce the \emph{relative superdifferential}
			\begin{align}\label{dfn:relative_superdifferential}
				\partial_{\textnormal{rel}}U(x)
					:=
					\begin{cases}
						[xU'_+(x), xU'_-(x)], & \textnormal{if }x>0,\\
						\{U'_+(0)\}, & \textnormal{if }x=0.
					\end{cases}
			\end{align}
			Note that this means the lower and upper endpoints of $\partial_{\textnormal{rel}}U(x)$ are measurable functions of $x$.
			
			For $x > 0$, this relative superdifferential simply equals $x \partial U(x)$.
			At the boundary $x = 0$, relative changes are not defined and every feasible direction is non-negative, so only the one-sided boundary slope $U'_+(0)$ will be relevant.
			We do not define the relative superdifferential at $\infty$, since every valid e-value is finite $Q_a$-almost surely and we only apply it to valid e-values.
			
			For a valid e-value $\varepsilon \in \mathcal{E}(H)$, we define its direction relative to $\varepsilon^*$ by
			\begin{align*}
				d_\varepsilon
					:=
				\begin{cases}
					\varepsilon,\quad & \textnormal{if } \varepsilon^*=0, \\
					\dfrac{\varepsilon-\varepsilon^*}{\varepsilon^*},\quad & \textnormal{if }\varepsilon^*>0.
				\end{cases}
			\end{align*}
			
			In the proof, we use the following standard compactness fact.
			
			\begin{lem}\label{lem:order_interval}
				Let $\underline{r}, \overline{r} \in L^1(Q_a)$ satisfy $0 \leq \underline{r} \leq \overline{r}$.
				Then the order interval
				\begin{align*}
					[\underline{r}, \overline{r}]
						:= \{r \in L^1(Q_a) : \underline{r} \leq r \leq \overline{r}\}
				\end{align*}
				is convex and weakly compact in $L^1(Q_a)$.
			\end{lem}
			\begin{proof}
				The order interval is uniformly integrable because all its elements are dominated by the integrable random variable $\overline{r}$.
				The Dunford--Pettis compactness criterion \citep{dunford1940linear} therefore implies that it is relatively weakly compact.
				Moreover, the order interval is convex and closed in $L^1(Q_a)$, and hence weakly closed.
				It is therefore weakly compact.
			\end{proof}
			
			\begin{proof}[Proof of Proposition \ref{prp:EU_converse_composite}]
				If $Q_a(\mathcal{X}) = 0$, the result is immediate with $\lambda^* = 0$ and $P^* = 0$.
				Suppose therefore that $Q_a(\mathcal{X}) > 0$.
				
				As explained in the main text before Proposition \ref{prp:EU_converse_composite}, concavity implies that the first-order condition already makes $\varepsilon^*$ $U$-optimal.
				It therefore remains to construct a single supporting measure $P^* \in H^\textnormal{eff}$ and establish the pointwise representation \eqref{opt:composite_U-optimal}.
				
				The proof has three parts.
				We first identify the compact set of relative supergradients over which the minimax argument will operate.
				We then use Sion's minimax theorem to select one relative supergradient that supports every feasible relative direction simultaneously.
				Finally, we convert this relative supergradient into an ordinary pointwise supergradient and normalize it into the supporting measure $P^*$. \\
								
				\noindent\textit{The compact set of relative supergradients}. 
					We start by defining the lower and upper endpoints of the relative superdifferential at $\varepsilon^*$ by
					\begin{align*}
						\underline{r}
							:= \inf \partial_{\textnormal{rel}}U(\varepsilon^*), \quad
						\overline{r}
							:= \sup \partial_{\textnormal{rel}}U(\varepsilon^*).
					\end{align*}
					These random variables are measurable by the established measurability of the endpoints of $\partial_{\textnormal{rel}}U(x)$.
					
					The central object in the first part is the set of all admissible relative supergradients
					\begin{align*}
						\mathcal{R}
							:= [\underline{r}, \overline{r}]
							= \{r \in L^1(Q_a) : r \in \partial_{\textnormal{rel}}U(\varepsilon^*),\ Q_a\textnormal{-a.s.}\}.
					\end{align*}
					This is the collection from which Sion's theorem will later select one common supporting element.
					Pointwise, $\partial_{\textnormal{rel}} U(\varepsilon^*)$ is a non-negative interval.
					Hence, to be able to apply Lemma \ref{lem:order_interval}, it suffices to verify that its upper endpoint $\overline{r}$ is $Q_a$-integrable.
					
					For this purpose, we apply the first-order condition to the e-value
					\begin{align*}
						1\{\varepsilon^* = 0\},
					\end{align*}
					which is valid because it is bounded from above by the constant e-value $1$.
					The corresponding direction from $\varepsilon^*$ is
					\begin{align*}
						1\{\varepsilon^*=0\} - \varepsilon^*
							=
						\begin{cases}
							1,
							&\textnormal{if }\varepsilon^*=0,\\
							-\varepsilon^*,
							&\textnormal{if }\varepsilon^*>0.
						\end{cases}
					\end{align*}
					Hence,
					\begin{align}\label{eq:special_direction_envelope}
						\mathrm{D}U\left(\varepsilon^*; 1\{\varepsilon^* = 0\} - \varepsilon^*\right)
							=
							\begin{cases}
								U'_+(0), &\textnormal{if }\varepsilon^*=0,\\
								-\varepsilon^*U'_-(\varepsilon^*), &\textnormal{if }\varepsilon^*>0.
							\end{cases}
					\end{align}
					As a consequence,
					\begin{align*}
						\left|\mathrm{D}U\left(\varepsilon^*; 1\{\varepsilon^*=0\} - \varepsilon^*\right)\right|
							= \overline{r},
					\end{align*}
					$Q_a$-almost surely.
					The directional derivative on the left is $Q_a$-integrable by assumption, so $\overline{r}$ is $Q_a$-integrable.
					Since $0 \leq \underline{r} \leq \overline{r}$, the lower endpoint $\underline{r}$ is also $Q_a$-integrable.
					Lemma \ref{lem:order_interval} therefore shows that $\mathcal{R}$ is non-empty, convex and weakly compact. \\

				\noindent\textit{Selecting one supergradient for all feasible directions}.
					The central object in the second part is the set of bounded feasible directions
					\begin{align*}
						\mathcal{D}
							:= \{d_\varepsilon : \varepsilon \in \mathcal{E}(H), d_\varepsilon \in L^\infty(Q_a)\}.
					\end{align*}
					We initially restrict to this bounded class so that it pairs continuously with the weakly compact set $\mathcal{R}$.
					The map $\varepsilon \mapsto d_\varepsilon$ is affine, so $\mathcal{D}$ is convex.
					It also contains the zero direction, corresponding to $\varepsilon = \varepsilon^*$.
					
					Fix a relative direction $d_\varepsilon \in \mathcal{D}$.
					The pointwise directional derivative of a concave function satisfies
					\begin{align}
						\mathrm{D}U(\varepsilon^*; \varepsilon - \varepsilon^*)
							&= \underline{r} d_\varepsilon 1\{d_\varepsilon \geq 0\} + \overline{r} d_\varepsilon 1\{d_\varepsilon < 0\} \nonumber\\
							&= \min_{r \in \partial_{\textnormal{rel}}U(\varepsilon^*)} r d_\varepsilon,\label{min:direction_dependent}
					\end{align}
					$Q_a$-almost surely.
					Indeed, a non-negative relative direction selects the lower endpoint $\underline{r}$, while a negative relative direction selects the upper endpoint $\overline{r}$.
					At zero, every feasible direction is non-negative and the two endpoints coincide with the boundary slope.
										
					Selecting $\underline{r}$ where $d_\varepsilon \geq 0$ and $\overline{r}$ where $d_\varepsilon < 0$ gives a measurable element of $\mathcal{R}$ attaining the pointwise minimum.
					Therefore, the first-order condition implies
					\begin{align}\label{ineq:FOC_rd}
						\min_{r \in \mathcal{R}} \textnormal{E}^{Q_a}[rd]
							= \textnormal{E}^{Q_a}[\mathrm{D}U(\varepsilon^*; \varepsilon - \varepsilon^*)]	
							\leq 0,
					\end{align}
					for $d = d_\varepsilon \in \mathcal{D}$.
					This means that, for every feasible relative direction $d \in \mathcal D$, there is some admissible relative supergradient $r \in \mathcal{R}$ such that
					\begin{align*}
						\textnormal{E}^{Q_a}[rd]
							\leq 0.
					\end{align*}
					This relative supergradient may depend on the direction $d$.
					To construct the measure $P^*$, we must show that one common element of $\mathcal{R}$ satisfies this inequality for every direction in $\mathcal{D}$ simultaneously.
					
					We obtain this common element by applying the minimax theorem of \citet{sion1958general} to the bilinear pairing
					\begin{align*}
						(r, d) \mapsto \textnormal{E}^{Q_a}[rd],
					\end{align*}
					where $\mathcal{R}$ carries the weak topology of $L^1(Q_a)$ and $\mathcal{D}$ the norm topology of $L^\infty(Q_a)$.
					The set $\mathcal{R}$ is compact and convex, while $\mathcal{D}$ is convex.
					Hence,
					\begin{align}\label{eq:minimax}
						\inf_{r \in \mathcal{R}} \sup_{d \in \mathcal{D}} \textnormal{E}^{Q_a}[rd]
							= \sup_{d \in \mathcal{D}} \min_{r \in \mathcal{R}} \textnormal{E}^{Q_a}[rd]
							\leq 0,
					\end{align}
					where the inequality follows from \eqref{ineq:FOC_rd}.
					On the other hand, since the zero direction belongs to $\mathcal{D}$,
					\begin{align*}
						\sup_{d \in \mathcal{D}} \textnormal{E}^{Q_a}[rd] \geq 0,
					\end{align*}
					for every $r \in \mathcal{R}$.
					As a consequence, both sides of \eqref{eq:minimax} equal zero.
					
					The map
					\begin{align*}
						r \mapsto \sup_{d \in \mathcal{D}} \textnormal{E}^{Q_a}[rd]
					\end{align*}
					is weakly lower semicontinuous, so it attains its minimum on the weakly compact set $\mathcal{R}$.
					We may therefore choose one common relative supergradient $r^* \in \mathcal{R}$ such that
					\begin{align}\label{ineq:one_gradient_to_rule_them_all}
						\textnormal{E}^{Q_a}[r^*d]
							\leq 0,
					\end{align}
					for every $d \in \mathcal{D}$.
					In other words, Sion's theorem has exchanged the direction-dependent pointwise choices in \eqref{min:direction_dependent} for one relative supergradient that supports every bounded feasible relative direction. \\
					
				\noindent\textit{Constructing the supporting measure}.
					We now convert the common relative supergradient back into an ordinary pointwise supergradient:
					\begin{align*}
						\gamma^*
							:=
							\begin{cases}
								r^*, \quad &\textnormal{ if } \varepsilon^* = 0, \\
								\frac{r^*}{\varepsilon^*}, \quad &\textnormal{ if } \varepsilon^* > 0.
							\end{cases}
					\end{align*}
					Since $r^* \in \mathcal{R}$, we have
					\begin{align*}
						r^*
							\in
						\partial_{\textnormal{rel}}U(\varepsilon^*),
					\end{align*}
					$Q_a$-almost surely.
					On $\{\varepsilon^*>0\}$, the definition of the relative superdifferential gives
					\begin{align*}
						r^*
							\in
						\varepsilon^*\partial U(\varepsilon^*).
					\end{align*}
					Dividing by $\varepsilon^*$ therefore shows that
					\begin{align}\label{inclusion:gamma_in_supergradient}
						\gamma^*
							= \frac{r^*}{\varepsilon^*}
							\in \partial U(\varepsilon^*).
					\end{align}
					On $\{\varepsilon^* = 0\}$, the relative superdifferential is the singleton $\{U'_+(0)\}$, so that
					\begin{align*}
						\gamma^* = r^*
							= U'_+(0).
					\end{align*}
					If this event $\{\varepsilon^* = 0\}$ has positive $Q_a$-measure, then $U(0)$ is finite because $U(\varepsilon^*)$ is $Q_a$-almost surely finite.
					Moreover, \eqref{eq:special_direction_envelope} and the assumed integrability of the corresponding directional derivative imply that $U'_+(0)<\infty$.
					Concavity then gives
					\begin{align*}
						U(z)-U(0)
							\leq U'_+(0)z, \quad z \geq 0,
					\end{align*}
					and hence $U'_+(0)\in\partial U(0)$.
					If $\{\varepsilon^*=0\}$ has zero $Q_a$-measure, its behavior is irrelevant.
					We have therefore shown that
					\begin{align}\label{ineq:common_ordinary_supergradient}
						\gamma^*
							\in \partial U(\varepsilon^*),
					\end{align}
					$Q_a$-almost surely.
					
					By the definitions of the relative direction and the descaled supergradient,
					\begin{align*}
						r^* d_\varepsilon
							= \gamma^* (\varepsilon - \varepsilon^*).
					\end{align*}
					As a consequence, \eqref{ineq:one_gradient_to_rule_them_all} gives
					\begin{align}\label{ineq:foc_bounded}
						\textnormal{E}^{Q_a}[\gamma^*(\varepsilon - \varepsilon^*)]
							\leq 0,
					\end{align}
					for every valid e-value with bounded relative direction.
					
					It remains to remove the boundedness restriction.
					Fix an arbitrary valid e-value $\varepsilon \in \mathcal{E}(H)$ and define the truncated e-value
					\begin{align*}
						\varepsilon_n
							:=
							\begin{cases}
								\varepsilon \wedge n, \quad &\textnormal{ if } \varepsilon^* = 0, \\
								\varepsilon \wedge n \varepsilon^*, \quad &\textnormal{ if } \varepsilon^* > 0.
							\end{cases}
					\end{align*}
					Since $\varepsilon_n \leq \varepsilon$, the random variable $\varepsilon_n$ remains a valid e-value.
					Its relative direction is bounded between $-1$ and $n$, so \eqref{ineq:foc_bounded} applies:
					\begin{align*}
						\textnormal{E}^{Q_a}[\gamma^* \varepsilon_n]
							\leq \textnormal{E}^{Q_a}[\gamma^*\varepsilon^*].
					\end{align*}
					Since $\varepsilon_n \uparrow \varepsilon$ and $\gamma^* \geq 0$, monotone convergence gives
					\begin{align}\label{dfn:lambda_star_proof}
						\textnormal{E}^{Q_a}[\gamma^*\varepsilon]
							\leq \lambda^*
							:=  \textnormal{E}^{Q_a}[\gamma^*\varepsilon^*]
					\end{align}
					for every $\varepsilon \in \mathcal{E}(H)$.
					The multiplier $\lambda^*$ is finite because
					\begin{align*}
						\gamma^*\varepsilon^*
							= r^* 1\{\varepsilon^* > 0\}
							\leq r^*
							\in L^1(Q_a).
					\end{align*}
					Taking the constant e-value $\varepsilon \equiv 1$ in \eqref{dfn:lambda_star_proof} also gives
					\begin{align*}
						\textnormal{E}^{Q_a}[\gamma^*]
							\leq \lambda^*.
					\end{align*}
					
					Suppose first that $\lambda^* > 0$.
					We then normalize the common pointwise supergradient into the non-negative measure $P^*$, defined by
					\begin{align}\label{dfn:P_star_from_supergradient}
						\frac{dP^*}{dQ_a}
							:= \frac{\gamma^*}{\lambda^*}.
					\end{align}
					The supporting inequality \eqref{dfn:lambda_star_proof} then gives
					\begin{align*}
						\textnormal{E}^{P^*}[\varepsilon]
							\leq 1,
					\end{align*}
					for every $\varepsilon \in \mathcal{E}(H)$.
					Hence, $P^* \in H^{\textnormal{eff}}$.
					Moreover, $\textnormal{E}^{P^*}[\varepsilon^*] = 1$.
					
					Combining \eqref{dfn:P_star_from_supergradient} and \eqref{ineq:common_ordinary_supergradient} gives
					\begin{align*}
						\lambda^* \frac{dP^*}{dQ_a}
							\in \partial U(\varepsilon^*),
					\end{align*}
					$Q_a$-almost surely.
					Equivalently, the e-value $\varepsilon^*$ satisfies
					\begin{align*}
						\varepsilon^*(x) \in \argmax_{z \in [0, \infty)} \left\{U(z) - \lambda^* \frac{dP^*}{dQ_a}(x) z\right\},
					\end{align*}
					$Q_a$-almost surely.
					The identity $\textnormal{E}^{P^*}[\varepsilon^*] = 1$ gives complementary slackness.
					
					Finally, suppose that $\lambda^* = 0$.
					Taking the constant e-value $\varepsilon \equiv 1$ in \eqref{dfn:lambda_star_proof} gives
					\begin{align*}
						\textnormal{E}^{Q_a}[\gamma^*]
							= 0,
					\end{align*}
					and hence $\gamma^* = 0$, $Q_a$-almost surely.
					By \eqref{ineq:common_ordinary_supergradient},
					\begin{align*}
						0 \in \partial U(\varepsilon^*),
					\end{align*}
					so the e-value $\varepsilon^*$ maximizes $U$ pointwise.
					Taking $P^* = 0$ proves the result also in this case.
			\end{proof}

		\subsection{Proof of Corollary \ref{cor:EU_converse_composite_radial}}
			\begin{proof}
				By Lemma \ref{lem:composite_well_defined}, $U(\varepsilon^*)$ is finite $Q_a$-almost surely.
				We show that the integrability condition permits the expected directional derivative to be obtained from $U$-optimality.
				
				Fix a valid e-value $\varepsilon \in \mathcal{E}(H)$ and define the line segment from $\varepsilon^*$ to this e-value by
				\begin{align*}
					\varepsilon_t
						:= (1 - t) \varepsilon^* + t \varepsilon,
				\end{align*}
				where $t \in (0, 1)$.
				Since $\mathcal{E}(H)$ is convex, every e-value $\varepsilon_t$ is valid.
				Optimality of $\varepsilon^*$ therefore gives
				\begin{align}\label{ineq:radial_directional_quotient}
					\textnormal{E}^{Q_a}\left[\frac{U(\varepsilon_t) - U(\varepsilon^*)}{t}\right]
						\leq 0.
				\end{align}
				
				To pass to the directional derivative as $t \downarrow 0$, we construct a common integrable lower bound for these difference quotients.
				For $0 < t < 1 - a$, define
				\begin{align*}
					\theta
						:= \frac{t}{1 - a}.
				\end{align*}
				Then
				\begin{align*}
					\varepsilon_t
						= (1 - \theta) \varepsilon^* + \theta (a \varepsilon^* + (1 - a) \varepsilon).
				\end{align*}
				Concavity and monotonicity of $U$ give
				\begin{align}\label{ineq:integrable_lower_bound}
					\frac{U(\varepsilon_t) - U(\varepsilon^*)}{t}
						\geq -\frac{U(\varepsilon^*) - U(a\varepsilon^*)}{1 - a}.
				\end{align}
				The right-hand side is $Q_a$-integrable by assumption.
				
				By concavity, the difference quotients in \eqref{ineq:radial_directional_quotient} increase pointwise as $t \downarrow 0$ to
				\begin{align*}
					\mathrm{D}U(\varepsilon^*; \varepsilon - \varepsilon^*).
				\end{align*}
				Since we have a common integrable lower bound \eqref{ineq:integrable_lower_bound}, the monotone convergence theorem gives
				\begin{align*}
					-\infty &< \textnormal{E}^{Q_a}\left[\mathrm{D}U(\varepsilon^*; \varepsilon - \varepsilon^*)\right] \\
						&= \lim_{t \downarrow 0}\textnormal{E}^{Q_a}\left[\frac{U(\varepsilon_t) - U(\varepsilon^*)}{t}\right]
						\leq 0.
				\end{align*}
				The same lower bound controls the negative part of the directional derivative, while the finiteness of the preceding limit controls its positive part.
				Hence, $\mathrm{D}U(\varepsilon^*; \varepsilon - \varepsilon^*)$ is $Q_a$-integrable.
				
				Since the valid e-value $\varepsilon$ was arbitrary, the integrability and first-order conditions of Proposition \ref{prp:EU_converse_composite} hold for every $\varepsilon \in \mathcal{E}(H)$.
				The proposition therefore yields the supporting measure $P^*$, the pointwise representation \eqref{opt:composite_U-optimal}, and the complementary slackness condition.
			\end{proof}
			
		\subsection{Proof of Theorem \ref{thm:existence_composite}}
			
			\begin{proof}[Proof of Theorem \ref{thm:existence_composite}]
				We prove the first claim: the second claim then follows from Lemma \ref{lem:composite_candidate_integrability}, given below.
				
				If $Q_a(\mathcal{X}) = 0$, then $\varepsilon_0$ is immediately $U$-optimal.
				Suppose therefore that $Q_a(\mathcal{X}) > 0$.
				
				We first describe the implications of the uniform integrability assumption.
				For every valid e-value $\varepsilon \in \mathcal{E}(H)$,
				\begin{align}\label{ineq:integrable_positive_part}
					\textnormal{E}^{Q_a}[[U(\varepsilon) - U(\varepsilon_0)]^+]	
						< \infty.
				\end{align}
				Indeed, if the positive part had infinite expectation, well-definedness of the relative utility objective would force
				\begin{align*}
					\textnormal{E}^{Q_a}[U(\varepsilon) - U(\varepsilon_0)]
						= \infty.
				\end{align*}
				This would imply $\varepsilon \in C(\varepsilon_0)$, contradicting uniform integrability.
				Moreover, for every $\varepsilon \in C(\varepsilon_0)$, the relative utility difference is $Q_a$-integrable, since its expected value is non-negative and its positive part is integrable.
				Since $\varepsilon_0 \in C(\varepsilon_0)$, uniform integrability therefore gives
				\begin{align*}
					0 
						&\leq \sup_{\varepsilon \in \mathcal{E}(H)} \textnormal{E}^{Q_a}[U(\varepsilon) - U(\varepsilon_0)] \\
						&= \sup_{\varepsilon \in C(\varepsilon_0)} \textnormal{E}^{Q_a}[U(\varepsilon) - U(\varepsilon_0)]
						< \infty.
				\end{align*}
				
				We now start with proving that an optimizer exists.
				For this purpose, we start by constructing a candidate optimizer from a maximizing sequence.
				Such a sequence need not itself converge, so we use the Koml\'os-type lemma to replace it by forward convex combinations that do converge.
				The convexity of $\mathcal{E}(H)$ ensures that these combinations remain valid, while Fatou's lemma will show that their limit is also valid.
				
				Choose a sequence $\varepsilon_n \in C(\varepsilon_0)$ such that
				\begin{align*}
					\textnormal{E}^{Q_a}[U(\varepsilon_n) - U(\varepsilon_0)]
						\to \sup_{\varepsilon \in \mathcal{E}(H)} \textnormal{E}^{Q_a}[U(\varepsilon) - U(\varepsilon_0)].
				\end{align*}
				Applying the Koml\'os-type Lemma 2.8 of \citet{larsson2024numeraire} under the probability measure $Q_a / Q_a(\mathcal{X})$, there exist forward convex combinations
				\begin{align*}
					\widetilde{\varepsilon}_n \in \textnormal{conv}(\varepsilon_n, \varepsilon_{n + 1}, \dots)
				\end{align*}
				that converge $Q_a$-almost surely.
				Define
				\begin{align*}
					\varepsilon^*
						:= \liminf_{n \to \infty} \widetilde{\varepsilon}_n
				\end{align*}
				everywhere.
				Since $\mathcal{E}(H)$ is convex, $\widetilde{\varepsilon}_n \in \mathcal{E}(H)$.
				Hence, for every $P \in H$, Fatou's lemma gives
				\begin{align*}
					\textnormal{E}^P[\varepsilon^*]
						\leq \liminf_{n \to \infty} \textnormal{E}^P[\widetilde{\varepsilon}_n]
						\leq 1.
				\end{align*}
				As a consequence, $\varepsilon^* \in \mathcal{E}(H)$.
				
				We next show that taking forward convex combinations has not destroyed the maximizing property.
				This will also keep the sequence inside $C(\varepsilon_0)$, where the uniform integrability assumption is available.
				
				Concavity of $U$ gives
				\begin{align*}
					\textnormal{E}^{Q_a}[U(\widetilde{\varepsilon}_n) - U(\varepsilon_0)]
						\geq \inf_{k \geq n} \textnormal{E}^{Q_a}[U(\varepsilon_k) - U(\varepsilon_0)]
						\geq 0.
				\end{align*}
				Hence, $\widetilde{\varepsilon}_n \in C(\varepsilon_0)$, and the forward convex combinations still form a maximizing sequence. In particular, the family
				\begin{align*}
					\left\{[U(\widetilde{\varepsilon}_n) - U(\varepsilon_0)]^+ : n \in \mathbb{N}\right\}
				\end{align*}
				is uniformly integrable under $Q_a$.
				
				We can now pass from the maximizing sequence to its limit.
				Upper semicontinuity preserves utility pointwise in the limit, while uniform integrability allows us to pass this inequality through the expectation.
				
				By upper semicontinuity of $U$,
				\begin{align*}
					U(\varepsilon^*) - U(\varepsilon_0)
						\geq \limsup_{n \to \infty} [U(\widetilde{\varepsilon}_n) - U(\varepsilon_0)],
				\end{align*}
				$Q_a$-almost surely.
				The reverse Fatou lemma, using uniform integrability of the positive parts, now gives
				\begin{align}
					\textnormal{E}^{Q_a}[U(\varepsilon^*) - U(\varepsilon_0)]
						&\geq \limsup_{n \to \infty} \textnormal{E}^{Q_a}[U(\widetilde{\varepsilon}_n) - U(\varepsilon_0)] \nonumber\\
						&= \sup_{\varepsilon \in \mathcal{E}(H)} \textnormal{E}^{Q_a}[U(\varepsilon) - U(\varepsilon_0)]. \label{ineq:supremum_optimal}
				\end{align}
				Since $\varepsilon^* \in \mathcal{E}(H)$, equality must hold.
				
				So far, we have shown optimality relative to the fixed benchmark $\varepsilon_0$.
				It remains to translate this benchmarked optimality into the full $U$-optimality.
				For this purpose, we first verify that the direct utility comparison between $\varepsilon^*$ and an arbitrary valid e-value is well-defined.
				
				Since the attained supremum in \eqref{ineq:supremum_optimal} is nonnegative, $\varepsilon^* \in C(\varepsilon_0)$, and therefore the difference $U(\varepsilon^*) - U(\varepsilon_0)$ is $Q_a$-integrable.
				For every $\varepsilon \in \mathcal{E}(H)$,
				\begin{align*}
					[U(\varepsilon^*) - U(\varepsilon)]^-
						\leq [U(\varepsilon) - U(\varepsilon_0)]^+
						+ [U(\varepsilon_0) - U(\varepsilon^*)]^+.
				\end{align*}
				The first term is integrable by \eqref{ineq:integrable_positive_part}, while the second is integrable because $U(\varepsilon^*) - U(\varepsilon_0)$ is integrable.
				Hence, the pairwise relative expected utility is well-defined.
				Moreover, the benchmark terms cancel, so that
				\begin{align*}
					\textnormal{E}^{Q_a}[U(\varepsilon^*) - U(\varepsilon)]
						&= \textnormal{E}^{Q_a}[U(\varepsilon^*) - U(\varepsilon_0)] \\
						&- \textnormal{E}^{Q_a}[U(\varepsilon) - U(\varepsilon_0)] \\
						&\geq 0.
				\end{align*}
				As a consequence, $\varepsilon^*$ is $U$-optimal.
			\end{proof}

			We separate the second claim of Theorem \ref{thm:existence_composite} into the following lemma.
			
			\begin{lem}\label{lem:composite_candidate_integrability}
				Suppose that $U$ is concave, upper semicontinuous and finite on $(0, \infty)$, and that
				\begin{align*}
					\limsup_{z \to \infty} \frac{U(z)}{z} \leq 0.
				\end{align*}
				Let $\varepsilon_0 \in \mathcal{E}(H)$ be so that $U(\varepsilon_0)$ is $Q_a$-integrable and define $C(\varepsilon_0)$ as in Theorem \ref{thm:existence_composite}.
				If $\sup_{\varepsilon \in C(\varepsilon_0)} \textnormal{E}^{Q_a}[\varepsilon] < \infty$, then
				\begin{align*}
					\left\{[U(\varepsilon) - U(\varepsilon_0)]^+ : \varepsilon \in C(\varepsilon_0)\right\}
				\end{align*}
				is uniformly integrable under $Q_a$.
			\end{lem}
			\begin{proof}
				The proof consists of two steps.
				We first combine the sublinear growth of $U$ with the assumed uniform bound on the $Q_a$-expectations of the candidate e-values to obtain uniform integrability of the positive utility gains relative to the constant benchmark $1$.
				We then transfer this conclusion to the chosen benchmark $\varepsilon_0$.

				The sublinear growth condition implies
				\begin{align}\label{lim:sublinear_difference}
					\frac{[U(z) - U(1)]^+}{z} \to 0,
				\end{align}
				as $z \to \infty$.
				Moreover, concavity and finiteness on $(0, \infty)$ rule out $U(0) = \infty$.
				Indeed, otherwise concavity would give
				\begin{align*}
					U(1) \geq \frac{1}{2} U(0) + \frac{1}{2}U(2) = \infty,
				\end{align*}
				contradicting finiteness of $U(1)$.
				Upper semicontinuity therefore implies that $z \mapsto [U(z) - U(1)]^+$ is bounded on compact intervals in $[0, \infty)$.
				
				We now use \eqref{lim:sublinear_difference} to control the utility tails.
				Fix $\eta > 0$ and choose $R$ sufficiently large that
				\begin{align*}
					[U(z) - U(1)]^+
						\leq \eta z,
				\end{align*}
				for every $z \geq R$.
				Define
				\begin{align*}
					K_R
						:= \sup_{0 \leq z \leq R} [U(z) - U(1)]^+
						< \infty.
				\end{align*}
				For every $K > K_R$,
				\begin{align*}
					\{[U(\varepsilon) - U(1)]^+ > K\}
						\subseteq \{\varepsilon > R\}.
				\end{align*}
				As a consequence,
				\begin{align*}
					\sup_{\varepsilon \in C(\varepsilon_0)} &\textnormal{E}^{Q_a}\left[[U(\varepsilon) - U(1)]^+ 1\{[U(\varepsilon) - U(1)]^+ > K\}\right] \\
						&\leq \eta \sup_{\varepsilon \in C(\varepsilon_0)} \textnormal{E}^{Q_a}[\varepsilon].
				\end{align*}
				Since the supremum on the right-hand side is finite and $\eta > 0$ is arbitrary, the family
				\begin{align*}
					\left\{[U(\varepsilon) - U(1)]^+ : \varepsilon \in C(\varepsilon_0)\right\}
				\end{align*}
				is uniformly integrable under $Q_a$.
				
				It now remains to transfer this conclusion from the constant benchmark $1$ to $\varepsilon_0$.
				For every $\varepsilon \in C(\varepsilon_0)$,
				\begin{align*}
					[U(\varepsilon) - U(\varepsilon_0)]^+
						\leq [U(\varepsilon) - U(1)]^+ + [U(1) - U(\varepsilon_0)]^+.
				\end{align*}
				The second term is one fixed $Q_a$-integrable random variable, since
				\begin{align*}
					[U(1) - U(\varepsilon_0)]^+
						\leq |U(1)| + |U(\varepsilon_0)|.
				\end{align*}
				Uniform integrability is preserved under addition of a fixed integrable random variable and under domination.
				This proves the result.
			\end{proof}

		\subsection{Proof of Corollary \ref{cor:log_growth_composite}}
			We first prove a technical lemma on the utility function, before proving the corollary itself.
			
			\begin{lem}\label{lem:non-decreasing_utility}
				If $U$ is concave, finite on $(0, \infty)$ and unbounded from above on $[0,\infty)$, then $U$ is non-decreasing on $[0, \infty)$.
			\end{lem}
			\begin{proof}
				For the sake of contradiction, suppose $U$ is unbounded from above on $[0,\infty)$ and $U(y) < U(x)$ for some $0 < x < y$.
				Then, concavity would imply, for all $z > y$,
				\begin{align*}
					\frac{U(z) - U(y)}{z - y}
						\leq \frac{U(y) - U(x)}{y - x}
						=:c
						< 0.
				\end{align*}
				Hence, $U(z) \leq U(y) + c (z - y)$.
				This implies $U$ is bounded from above on $[y, \infty)$.
				It remains to control $U$ near zero.
				For every $0 < t < x$, concavity similarly implies
				\begin{align*}
					\frac{U(x) - U(t)}{x - t}
						\geq \frac{U(y) - U(x)}{y - x}
						= c.
				\end{align*}
				Rearranging gives
				\begin{align*}
					U(t)
						\leq U(x) - c (x - t)
						\leq U(x) - cx,
				\end{align*}
				where the final inequality uses $c < 0$ and $x - t \leq x$.
				As a consequence, $U$ is bounded from above on $(0, x)$.
				Finally, since a finite concave function is continuous on $(0, \infty)$, it is bounded from above on the compact interval $[x, y]$.
				Moreover, $U(0) = \infty$ is impossible, since concavity would give, for every $x > 0$,
				\begin{align*}
				    U(x / 2)
				        \geq \frac{1}{2} U(0) + \frac{1}{2} U(x)
				        = \infty,
				\end{align*}
				contradicting the finiteness of $U$ on $(0, \infty)$.
				As a consequence, $U$ is bounded from above on $[0,\infty)$, which is a contradiction.
				Therefore, $U$ is non-decreasing on $(0,\infty)$.

				It remains to include the endpoint zero.
				If $U(0) = -\infty$, then $U(0) \leq U(x)$ for every $x > 0$.
				This means we only still need to handle the situation in which $U(0)$ is finite.
				Suppose, for the sake of contradiction, that $U(0) > U(x)$ for some $x > 0$, and define
				\begin{align*}
    				c  := \frac{U(x)-U(0)}{x}  < 0.
				\end{align*}
				Then concavity implies, for every $z > x$, that
				\begin{align*}
				    \frac{U(z) - U(x)}{z-x}
				        \leq c,
				\end{align*}
				and hence
				\begin{align*}
				    U(z)
				        \leq U(x)+c(z-x)
				        \leq U(x).
				\end{align*}
				Thus, $U$ is bounded from above on $[x, \infty)$.
				Since $U$ is already known to be non-decreasing on $(0,\infty)$, it is also bounded from above by $U(x)$ on $(0,x]$.
				Together with the finiteness of $U(0)$, this implies that $U$ is bounded from above on $[0,\infty)$, which is a contradiction.
				Therefore, $U(0)\leq U(x)$ for every $x>0$, and $U$ is non-decreasing on $[0, \infty)$.	
			\end{proof}

			We are now ready to prove the corollary.
			
			\begin{proof}[Proof of Corollary \ref{cor:log_growth_composite}]
				If $Q_a(\mathcal{X}) = 0$, the constant e-value $1$ is $U$-optimal.
				Suppose therefore that $Q_a(\mathcal{X}) > 0$ and define $\overline{Q}_a := Q_a / Q_a(\mathcal{X})$.
				This rescaling does not affect $U$-optimality or uniform integrability.
				
				We first consider the case where $U$ is bounded from above.
				Let
				\begin{align*}
					M := \sup_{z \in [0, \infty)} U(z) < \infty,
				\end{align*}
				and consider the benchmark e-value $\varepsilon_0 \equiv 1$.
				Then $\varepsilon_0 \in \mathcal{E}(H)$ and $U(1)$ is finite.
				Moreover, for every $\varepsilon \in \mathcal{E}(H)$, we have the uniform bound
				\begin{align*}
					[U(\varepsilon) - U(1)]^+
						\leq \sup_{z \geq 0} U(z) - U(1)
						= M - U(1).
				\end{align*}
				It follows both that
				\begin{align*}
					\textnormal{E}^{Q_a}[U(\varepsilon) - U(1)]
				\end{align*}
				is well-defined as an element of $[-\infty, \infty)$ for every $\varepsilon \in \mathcal{E}(H)$ and that
				\begin{align*}
					\left\{[U(\varepsilon) - U(1)]^+ : \varepsilon \in C(1)\right\}
				\end{align*}
				is uniformly integrable under $Q_a$.
				As a consequence, the first part of Theorem \ref{thm:existence_composite} applies.

				We may therefore assume that $U$ is not bounded from above.
				By Lemma \ref{lem:non-decreasing_utility}, $U$ must then be non-decreasing on $[0, \infty)$.
								
				By the construction of the generalized Lebesgue decomposition, $\overline{Q}_a$ assigns zero mass to every $H$-negligible event.
				Hence, by Theorem 2.6 in \citet{larsson2024numeraire}, there exists a num\'eraire e-value $\varepsilon^\textnormal{num} \in \mathcal{E}(H)$ satisfying $0 < \varepsilon^\textnormal{num} < \infty$ $\overline{Q}_a$-almost surely, and
				\begin{align}
					\textnormal{E}^{\overline{Q}_a}\left[\frac{\varepsilon}{\varepsilon^\textnormal{num}}\right]
						\leq 1,
				\end{align}
				for every $\varepsilon \in \mathcal{E}(H)$.
				
				Define, only for the sake of this proof, the benchmark e-value
				\begin{align*}
					\varepsilon_0
						:= \frac{1 + \varepsilon^\textnormal{num}}{2}.
				\end{align*}
				Then $\varepsilon_0 \in \mathcal{E}(H)$, since both 1 and $\varepsilon^\textnormal{num}$ belong to $\mathcal{E}(H)$ and this class is convex.
				Moreover,
				\begin{align*}
					\frac{1}{2} \leq \varepsilon_0 < \infty,
				\end{align*}
				$\overline{Q}_a$-almost surely so that $U(\varepsilon_0)$ is $\overline{Q}_a$-almost surely finite.
				Since $\varepsilon_0 \geq \varepsilon^\textnormal{num}/2$,
				\begin{align}\label{ineq:relative_first_moment}
					\sup_{\varepsilon \in \mathcal{E}(H)} \textnormal{E}^{\overline Q_a} \left[\frac{\varepsilon}{\varepsilon_0}\right]
    					&\leq 2 \sup_{\varepsilon \in \mathcal{E}(H)}\textnormal{E}^{\overline Q_a}\left[\frac{\varepsilon}{\varepsilon^{\mathrm{num}}}\right]
    					\leq2.
				\end{align}
				
				Define
				\begin{align*}
					A := U(1) - U(1/2) < \infty.
				\end{align*}
				We claim that
				\begin{align}\label{ineq:extended_log_growth}
					[U(y) - U(x)]^+
						\leq A + K \log(y/x),
				\end{align}
				for all $1/2 \leq x \leq y < \infty$.
				Indeed, if $y \leq 1$, monotonicity gives
				\begin{align*}
					0 
						\leq U(y) - U(x)
						\leq U(1) - U(1/2)
						= A.
				\end{align*}
				If $x \geq 1$, then this follows directly from \eqref{ineq:log_growth}.
				Finally, if $x < 1 < y$ then
				\begin{align*}
					U(y) - U(x)
						&= [U(y) - U(1)] + [U(1) - U(x)] \\
						&\leq K \log y + A \\
						&\leq K \log(y/x) + A,
				\end{align*}
				where the first inequality uses \eqref{ineq:log_growth} to bound the first term and the definition of $A$ to bound the second term, and the second inequality uses $x \leq 1$.
				This proves \eqref{ineq:extended_log_growth}.
				
				Since $U$ is non-decreasing, the difference $U(\varepsilon) - U(\varepsilon_0)$ can be positive only when $\varepsilon > \varepsilon_0$.
				Applying \eqref{ineq:extended_log_growth} on this event gives
				\begin{align}\label{ineq:utility_difference_bound}
					[U(\varepsilon) - U(\varepsilon_0)]^+
						\leq A + K[\log(\varepsilon / \varepsilon_0)]^+.
				\end{align}
				
				By \eqref{ineq:relative_first_moment},
				\begin{align*}
					\sup_{\varepsilon \in \mathcal{E}(H)} \textnormal{E}^{\overline{Q}_a}\left[\frac{\varepsilon}{\varepsilon_0}\right] \leq 2.
				\end{align*}
				We claim that
				\begin{align*}
					\left\{\left[\log\left(\frac{\varepsilon}{\varepsilon_0}\right)\right]^+ : \varepsilon \in \mathcal{E}(H)\right\}
				\end{align*}
				is uniformly integrable.
				Indeed, fix $R \geq 1$.
				The function $z \mapsto (\log(z)) / z$ is decreasing on $[\exp(R), \infty)$.
				Therefore,
				\begin{align*}
					\frac{\log\left(\varepsilon / \varepsilon_0\right)}{\varepsilon / \varepsilon_0}
						\leq \frac{\log(\exp(R))}{\exp(R)}
						= R \exp(-R),
				\end{align*}
				on $[\exp(R), \infty)$.
				As a consequence,
				\begin{align*}
					\sup_{\varepsilon \in \mathcal{E}(H)} &\textnormal{E}^{\overline{Q}_a}\left[\left\{\log\left(\frac{\varepsilon}{\varepsilon_0}\right)\right\}^+1\left\{\left[\log\left(\frac{\varepsilon}{\varepsilon_0}\right)\right]^+ > R\right\}\right] \\
						&\leq R \exp(-R) \sup_{\varepsilon \in \mathcal{E}(H)} \textnormal{E}^{\overline{Q}_a}\left[\frac{\varepsilon}{\varepsilon_0}\right] \\
						&\leq 2 R \exp(-R).
				\end{align*}
				Since $2 R \exp(-R) \to 0$ as $R \to \infty$, the uniform integrability follows.
				
				By \eqref{ineq:utility_difference_bound}, this means that
				\begin{align*}
					\left\{\left[U(\varepsilon) - U(\varepsilon_0)\right]^+ : \varepsilon \in \mathcal{E}(H)\right\}
				\end{align*}
				is also uniformly integrable.
				In particular, the relative objective
				\begin{align*}
					\textnormal{E}^{Q_a}\left[U(\varepsilon) - U(\varepsilon_0)\right],
				\end{align*}
				is well-defined as an element of $[-\infty, \infty)$, for every $\varepsilon \in \mathcal{E}(H)$.
				Moreover, the required uniform integrability holds in the smaller subfamily $C(\varepsilon_0)$.
				The first part of Theorem \ref{thm:existence_composite} then gives the existence of the $U$-optimal e-value.
			\end{proof}

		\subsection{Proof of Proposition \ref{prp:composite_optimizer_properties}}
			\begin{proof}
				The uniqueness proof uses only strict concavity of $U$ and convexity of the feasible class, and $\mathcal{E}(H)$ is convex. 
				The positivity conclusion follows from Lemma \ref{lem:composite_well_defined}.
				
				For strict monotonicity, suppose that
				\begin{align*}
					\sup_{P\in H} \textnormal{E}^P[\varepsilon^*] < 1.
				\end{align*}
				Then
				\begin{align*}
					\varepsilon^*+\frac{1 - \sup_{P\in H} \textnormal{E}^P[\varepsilon^*]}{2}
				\end{align*}
				is a valid e-value and has strictly larger utility $Q_a$-almost surely, contradicting $U$-optimality.
			\end{proof}
		
		\subsection{Proof of Corollary \ref{cor:composite_structure}}
			\begin{proof}
				The pointwise condition in Theorem \ref{thm:EU_NP_composite} is the same scalar optimization problem as in the simple setting after the stated substitution.
				All conclusions therefore follow from the corresponding arguments in Corollary \ref{cor:structure_pointwise}.
				For strictly increasing $U$, complementary slackness gives $\textnormal{E}^{P^*}[\varepsilon^*] = 1$ once $\lambda^*>0$.
			\end{proof}

		\subsection{Proof of Lemma \ref{lem:composite_well_defined}}
			\begin{proof}
				The first implication in Lemma \ref{lem:well-defined} is unchanged.
				Indeed, taking $\varepsilon = \varepsilon^*$ shows that
				\begin{align*}
					U(\varepsilon^*) - U(\varepsilon^*)
				\end{align*}	
				can be well-defined only if $U(\varepsilon^*)$ is finite $Q_a$-almost surely.
				
				For the converse, fix $\varepsilon \in \mathcal{E}(H)$.
				Pointwise optimality gives
				\begin{align*}
					U(\varepsilon) - U(\varepsilon^*)
						\leq \lambda^* \frac{dP^*}{dQ_a} (\varepsilon - \varepsilon^*), \quad Q_a\textnormal{-a.s.}
				\end{align*}
				As a consequence,
				\begin{align*}
					\left[U(\varepsilon) - U(\varepsilon^*)\right]^+
						\leq \lambda^* \frac{dP^*}{dQ_a} \varepsilon.
				\end{align*}
				Since $P^* \in H^{\textnormal{eff}}$,
				\begin{align*}
					\textnormal{E}^{Q_a}\left[\frac{dP^*}{dQ_a}\varepsilon\right]
						= \textnormal{E}^{P^*}[\varepsilon]
						\leq 1.
				\end{align*}
				Hence, the positive part is integrable and the relative expected utility is well-defined.
				
				It remains to verify the composite version of Lemma \ref{lem:well-defined_U}.
				Since $\varepsilon^*$ is valid, the event $\{\varepsilon^* = \infty\}$ is $H$-negligible, and therefore $\varepsilon^* < \infty$, $Q_a$-almost surely.
				If $U(z) < \infty$ for every finite $z$, it follows that $U(\varepsilon^*) < \infty$, $Q_a$-almost surely.
				
				Let $z_0 < \infty$ satisfy $U(z_0) > -\infty$.
				Pointwise optimality gives
				\begin{align*}
					U(\varepsilon^*) - \lambda^* \frac{dP^*}{dQ_a} \varepsilon^*
						\geq U(z_0) - \lambda^* \frac{dP^*}{dQ_a} z_0
						> -\infty,
				\end{align*}
				$Q_a$-almost surely.
				Since $\varepsilon^*$ and $dP^*/dQ_a$ are finite $Q_a$-almost surely, this implies $U(\varepsilon^*) > -\infty$ $Q_a$-almost surely.
				As a consequence, $U(\varepsilon^*)$ is finite $Q_a$-almost surely, and the preceding part proves well-definedness.
			\end{proof}
			
		\subsection{Proof of Corollary \ref{cor:NP}}
			\begin{proof}
				The utility $U_\alpha^{\textnormal{NP}}$ is concave, upper semicontinuous and finite on $(0,\infty)$, and satisfies the logarithmic growth condition.
				Corollary \ref{cor:log_growth_composite} therefore gives the existence of an optimizer.
				Moreover, replacing any valid e-value $\varepsilon$ by $\varepsilon \wedge 1/\alpha$
				preserves its utility and validity, so an optimizer can be chosen to take values in $[0,1/\alpha]$.
				
				For the equivalence, suppose first that $\varepsilon^{\textnormal{NP}}$ is optimal.
				For every $a \in(0,1)$,
				\begin{align*}
					U_\alpha^{\textnormal{NP}}(\varepsilon^{\textnormal{NP}}) - U_\alpha^{\textnormal{NP}}(a\varepsilon^{\textnormal{NP}})
						= (1 - a) \varepsilon^{\textnormal{NP}}
						\leq \frac{1 - a}{\alpha}.
				\end{align*}
				Hence, Corollary \ref{cor:EU_converse_composite_radial} gives $\lambda^*\in[0,\infty)$ and $P^*\in H^{\textnormal{eff}}$ with $P^*\ll Q_a$ such that complementary slackness holds and $\varepsilon^{\textnormal{NP}}$ pointwise maximizes
				\begin{align*}
					z \wedge 1 / \alpha
						- \lambda^* \frac{dP^*}{dQ_a}(x) z.
				\end{align*}
				
				For a fixed $x$, write
				\begin{align*}
					y := \lambda^* \frac{dP^*}{dQ_a}(x).
				\end{align*}
				On $[0, 1/\alpha]$, the pointwise objective equals $(1-y)z$, while beyond $1/\alpha$ it equals $1/\alpha - yz$.
				Its capped maximizers are therefore
				\begin{align*}
					\begin{cases}
						1/\alpha, & y < 1,\\
						[0,1/\alpha], & y = 1,\\
						0, & y > 1.
					\end{cases}
				\end{align*}
				This proves the second statement.
				
				Conversely, suppose the second statement holds.
				The preceding scalar calculation shows that the optimizer in the second statement is precisely the pointwise maximization condition in Theorem \ref{thm:EU_NP_composite}.
				Since $\varepsilon^{\textnormal{NP}}$ is valid, complementary slackness holds and $U_\alpha^{\textnormal{NP}} (\varepsilon^{\textnormal{NP}})$ is finite, that theorem implies that $\varepsilon^{\textnormal{NP}}$ is $U_\alpha^{\textnormal{NP}}$-optimal.
			\end{proof}
			
		\subsection{Proof of Corollary \ref{cor:log_utility}}
			\begin{proof}
				The existence follows from Corollary \ref{cor:log_growth_composite}, since log-utility satisfies the logarithmic growth condition.
				The uniqueness and positivity conclusions of Proposition \ref{prp:composite_optimizer_properties} give $Q_a$-almost sure strict positivity and uniqueness.
				
				It remains to prove the equivalence of the three statements.
				We start by proving that the first statement implies the second.
				
				Suppose that $\varepsilon^\textnormal{log}$ is $U^\textnormal{log}$-optimal.
				For every $a \in (0, 1)$,
				\begin{align*}
					U^\textnormal{log}(\varepsilon^\textnormal{log}) - U^\textnormal{log}(a\varepsilon^\textnormal{log})
						= -\log(a).
				\end{align*}
				Corollary \ref{cor:EU_converse_composite_radial} therefore gives the integrable first-order condition \eqref{ineq:integrable_directional_optimality}.
				For log utility,
				\begin{align}\label{eq:log_directional_derivative}
					\mathrm{D}U^\textnormal{log}(\varepsilon^\textnormal{log}; \varepsilon - \varepsilon^\textnormal{log})
						= \frac{\varepsilon}{\varepsilon^\textnormal{log}} - 1.
				\end{align}
				Hence,
				\begin{align*}
					\textnormal{E}^{Q_a}\left[\frac{\varepsilon}{\varepsilon^\textnormal{log}} - 1\right]
						\leq 0,
				\end{align*}
				for every $\varepsilon \in \mathcal{E}(H)$, proving the second statement.
				
				Now suppose that the second statement holds.
				Taking the constant e-value $\varepsilon\equiv1$ shows that
				\begin{align*}
					\varepsilon^{\textnormal{log}}
						>0,
				\end{align*}
				$Q_a$-almost surely, since otherwise $1/\varepsilon^{\textnormal{log}}=\infty$ on an event of positive $Q_a$-measure.
				This means that, for every $\varepsilon \in \mathcal{E}(H)$,
				\begin{align*}
					\frac{\varepsilon}{\varepsilon^{\textnormal{log}}} - 1
				\end{align*}
				is well-defined, $Q_a$-integrable and has non-positive expectation.
				The conditions of Proposition \ref{prp:EU_converse_composite} are therefore satisfied.
				The proposition gives some $\lambda^* \in [0, \infty)$ and $P^* \in H^\textnormal{eff}$ with $P^* \ll Q_a$ such that the pointwise condition and complementary slackness hold.
				
				The positive multiplier and supergradient form conclusions of Corollary \ref{cor:composite_structure} give $\lambda^* > 0$ and 
				\begin{align}\label{eq:supergradient_log}
					\frac{1}{\varepsilon^\textnormal{log}}
						= \lambda^* \frac{dP^*}{dQ_a},
				\end{align}
				$Q_a$-almost surely.
				Complementary slackness and $\lambda^* > 0$ give $\textnormal{E}^{P^*}[\varepsilon^\textnormal{log}] = 1$.
				
				Multiplying by $\varepsilon^\textnormal{log}$ and integrating with respect to $Q_a$ yields
				\begin{align*}
					Q_a(\mathcal{X})
						= \lambda^* \textnormal{E}^{P^*}[\varepsilon^\textnormal{log}]
						= \lambda^*.
				\end{align*}
				Hence,
				\begin{align}\label{eq:log_strict_positive_PQ}
					\frac{dP^*}{dQ_a}
						= \frac{1}{Q_a(\mathcal{X}) \varepsilon^\textnormal{log}}
						> 0,
				\end{align}
				$Q_a$-almost surely.
				Since $P^* \ll Q_a$, the strict positivity in \eqref{eq:log_strict_positive_PQ} also gives $Q_a \ll P^*$.
				Hence, $P^*$ and $Q_a$ are mutually absolutely continuous.
				
				We may therefore invert the Radon--Nikodym derivative in \eqref{eq:log_strict_positive_PQ} to obtain
				\begin{align*}
					\varepsilon^\textnormal{log}
						= \frac{1}{Q_a(\mathcal{X})} \frac{dQ_a}{dP^*},
				\end{align*}
				$Q_a$-almost surely.
				Together with $\textnormal{E}^{P^*}[\varepsilon^\textnormal{log}] = 1$, this proves the third statement.
				
				Finally, suppose that the third statement holds.
				Define
				\begin{align*}
					\lambda^*
						:= Q_a(\mathcal{X}).
				\end{align*}
				Since $P^*$ and $Q_a$ are mutually absolutely continuous, the likelihood ratio representation can be inverted to give
				\begin{align*}
					\lambda^* \frac{dP^*}{dQ_a}
						= \frac{1}{\varepsilon^\textnormal{log}},
				\end{align*}
				$Q_a$-almost surely. 
				For every $x$, the function
				\begin{align*}
					z \mapsto \log(z) - \lambda^* \frac{dP^*}{dQ_a}(x) z
				\end{align*}
				is strictly concave, and its derivative vanishes at $z = \varepsilon^{\textnormal{log}}$.
				As a consequence,
				\begin{align*}
					\varepsilon^\textnormal{log}(x)
						\in \argmax_{z \in [0, \infty)} \left\{\log(z) - \lambda^* \frac{dP^*}{dQ_a}(x) z\right\},
				\end{align*}
				$Q_a$-almost surely.
				Moreover, $\lambda^* (1 - \textnormal{E}^{P^*}[\varepsilon^\textnormal{log}]) = 0$.
				Hence, all conditions of Theorem \ref{thm:EU_NP_composite} are satisfied, so that $\varepsilon^\textnormal{log}$ is $U^\textnormal{log}$-optimal.
			\end{proof}
			
		\subsection{Proof of Corollary \ref{cor:capped_power_characterization}}
			\begin{proof}
				We show that the first statement implies both the second and third statements, and that each of the latter statements implies the first.
				
				Suppose first that $\varepsilon^{(\alpha,h)}$ is $U_{\alpha,h}$-optimal.
				We begin by verifying the radial integrability condition in Corollary \ref{cor:EU_converse_composite_radial}.
				
				If $h\leq0$, the Positivity conclusion of Proposition \ref{prp:composite_optimizer_properties} gives
				\begin{align*}
					\varepsilon^{(\alpha,h)}
						>0,
				\end{align*}
				$Q_a$-almost surely.
				For $h\neq0$, we then have, for every $a\in(0,1)$,
				\begin{align}\label{eq:capped_power_radial_loss}
					U_{\alpha,h}\left(\varepsilon^{(\alpha,h)}\right) - U_{\alpha,h}\left(a\varepsilon^{(\alpha,h)}\right)
						= \frac{1-a^h}{h} \left(\varepsilon^{(\alpha,h)}\right)^h.
				\end{align}
				For $h=0$, the corresponding difference equals
				\begin{align*}
					U_{\alpha,0}\left(\varepsilon^{(\alpha,0)}\right) - U_{\alpha,0}\left(a\varepsilon^{(\alpha,0)}\right)
						=-\log(a).
				\end{align*}
				
				If $h>0$, the right-hand side of \eqref{eq:capped_power_radial_loss} is integrable because $\varepsilon^{(\alpha,h)}$ is bounded when $\alpha>0$, while its integrability is assumed when $\alpha=0$.
				If $h<0$, optimality relative to the constant e-value $1$ gives
				\begin{align*}
					\textnormal{E}^{Q_a}\left[U_h\left(\varepsilon^{(\alpha,h)}\right)\right]
						\geq 0.
				\end{align*}
				Since $U_h$ is bounded from above, its positive part is integrable, while the non-negative expected value rules out an infinite negative part.
				Hence, $U_h\left(\varepsilon^{(\alpha,h)}\right)$ and therefore $\left(\varepsilon^{(\alpha,h)}\right)^h$ are $Q_a$-integrable.
				Thus, the radial integrability condition holds in every case.
				
				Corollary \ref{cor:EU_converse_composite_radial} therefore gives the integrable first-order condition, together with some $\lambda^*\in[0,\infty)$ and $P^*\in H^{\textnormal{eff}}$ with $P^*\ll Q_a$ satisfying the pointwise condition and complementary slackness.
				
				Applying the integrable first-order condition to the constant e-value $1$ shows that
				\begin{align*}
					\varepsilon^{(\alpha,h)}
						>0,
				\end{align*}
				$Q_a$-almost surely, since the directional derivative of $U_{\alpha,h}$ from zero towards one is infinite for every $h<1$.
				
				Now let $\varepsilon$ be a valid e-value taking values in $[0,1/\alpha]$.
				The line segment from $\varepsilon^{(\alpha,h)}$ to $\varepsilon$ remains in $[0,1/\alpha]$, where $U_{\alpha,h}$ agrees with $U_h$.
				Consequently,
				\begin{align*}
					\mathrm{D}U_{\alpha,h}\left(\varepsilon^{(\alpha,h)};\varepsilon-\varepsilon^{(\alpha,h)}\right)
						= \left(\varepsilon^{(\alpha,h)}\right)^{h-1}\left(\varepsilon-\varepsilon^{(\alpha,h)}\right).
				\end{align*}
				The integrable first-order condition therefore gives \eqref{ineq:capped_power_first_order}, which proves the second statement.
				
				To obtain the third statement, fix $x\in\mathcal X$ and write
				\begin{align*}
					y
						:= \lambda^*\frac{dP^*}{dQ_a}(x).
				\end{align*}
				On $(0,1/\alpha)$, the derivative of the pointwise objective
				\begin{align*}
					z \mapsto U_{\alpha,h}(z) - yz
				\end{align*}
				equals
				\begin{align*}
					z^{h-1} - y.
				\end{align*}
				Since $h<1$, this derivative is strictly decreasing, while the pointwise objective is non-increasing beyond the cap.
				If $\alpha = 0$, pointwise optimality forces $y > 0$, since for $y = 0$ the objective is strictly increasing and has no finite maximizer.
				Its maximizer among the capped values is therefore
				\begin{align*}
					y^{-1/(1-h)}\wedge\frac{1}{\alpha}.
				\end{align*}
				The pointwise condition supplied by Corollary \ref{cor:EU_converse_composite_radial} therefore gives \eqref{opt:capped_power_utility}.
				The same corollary gives
				\begin{align*}
					\lambda^*\left(1-\textnormal{E}^{P^*}[\varepsilon^{(\alpha,h)}]\right)
						=0,
				\end{align*}
				which proves the third statement.
				
				We next show that the second statement implies the first.
				Since $1\leq1/\alpha$, the constant e-value $1$ is among the e-values considered in the second statement.
				The required integrability for this e-value forces
				\begin{align*}
					\varepsilon^{(\alpha,h)}
						>0,
				\end{align*}
				$Q_a$-almost surely.
				
				Let $\varepsilon\in\mathcal E(H)$ be arbitrary and define
				\begin{align*}
					\overline{\varepsilon}
						:= \varepsilon\wedge1/\alpha.
				\end{align*}
				Then $\overline{\varepsilon}$ is a valid e-value taking values in $[0,1/\alpha]$, and
				\begin{align*}
					U_{\alpha,h}(\varepsilon)
						=
					U_h(\overline{\varepsilon}).
				\end{align*}
				Concavity of $U_h$ gives
				\begin{align*}
					U_{\alpha,h}(\varepsilon)-U_{\alpha,h}\left(\varepsilon^{(\alpha,h)}\right)
						&=
					U_h(\overline{\varepsilon})-U_h\left(\varepsilon^{(\alpha,h)}\right) \\
						&\leq
					\left(\varepsilon^{(\alpha,h)}\right)^{h-1}\left(\overline{\varepsilon}-\varepsilon^{(\alpha,h)}\right).
				\end{align*}
				The right-hand side is $Q_a$-integrable and has non-positive expectation by the second statement.
				It follows that the relative expected utility is well-defined and
				\begin{align*}
					\textnormal{E}^{Q_a}\left[U_{\alpha,h}(\varepsilon)-U_{\alpha,h}\left(\varepsilon^{(\alpha,h)}\right)\right]
						\leq 0.
				\end{align*}
				Since the valid e-value $\varepsilon$ was arbitrary, $\varepsilon^{(\alpha,h)}$ is $U_{\alpha,h}$-optimal.
				
				Finally, suppose that the third statement holds.
				The scalar pointwise calculation above shows that \eqref{opt:capped_power_utility} is precisely the pointwise maximization condition in Theorem \ref{thm:EU_NP_composite}.
				The representation and validity of $\varepsilon^{(\alpha,h)}$ imply
				\begin{align*}
					0
						< \varepsilon^{(\alpha,h)}
						< \infty,
				\end{align*}
				$Q_a$-almost surely, so that $U_{\alpha,h}\left(\varepsilon^{(\alpha,h)}\right)$ is finite $Q_a$-almost surely.
				Together with complementary slackness, all conditions of Theorem \ref{thm:EU_NP_composite} are satisfied.
				That theorem therefore implies that $\varepsilon^{(\alpha,h)}$ is $U_{\alpha,h}$-optimal.
			\end{proof}

		\subsection{Proof of Corollary \ref{cor:capped_power_existence}}
			\begin{proof}
				If $\alpha > 0$ or $h \leq 0$, the utility $U_{\alpha,h}$ satisfies the logarithmic growth condition.
				Indeed, below the cap,
				\begin{align*}
					zU_h'(z)
						= z^h,
				\end{align*}
				which is bounded on $[1, 1/\alpha]$ if $\alpha > 0$ and on $[1, \infty)$ if $h \leq 0$.
				Above the cap, the utility is constant.
				Existence therefore follows from Corollary \ref{cor:log_growth_composite}.
				If $\alpha > 0$, truncating an e-value at $1/\alpha$ preserves its utility and validity.
				
				Now suppose that $\alpha = 0$ and $0 < h < 1$.
				We apply Theorem \ref{thm:existence_composite} with benchmark $\varepsilon_0$.
				Since
				\begin{align*}
					U_h(\varepsilon) - U_h(\varepsilon_0)
						= \frac{\varepsilon^h - \varepsilon_0^h}{h},
				\end{align*}
				the candidate class in that theorem is precisely $C_h(\varepsilon_0)$.
				The assumption $\textnormal{E}^{Q_a}[\varepsilon_0^h] < \infty$ also ensures that $U_h(\varepsilon_0)$ is $Q_a$-integrable and that all relative utility comparisons with $\varepsilon_0$ are well-defined.
				Moreover, $\lim_{z \to \infty} U_h(z) / z = 0$.
				The second part of Theorem \ref{thm:existence_composite} therefore gives the second statement.
				
				Finally, suppose that $\alpha = 0$ and $h = 1$.
				Then $U_{0,1}(z) = z - 1$, and hence
				\begin{align*}
					U_{0,1}(\varepsilon) - U_{0,1}(\varepsilon_0)
						= \varepsilon - \varepsilon_0.
				\end{align*}
				Moreover, $(\varepsilon-\varepsilon_0)^+ \leq \varepsilon \leq (\varepsilon-\varepsilon_0)^+ + \varepsilon_0$.
				Since $\varepsilon_0$ is $Q_a$-integrable, the family
				\begin{align*}
					\left\{(\varepsilon - \varepsilon_0)^+: \varepsilon \in C_1(\varepsilon_0)\right\}
				\end{align*}
				is uniformly integrable if and only if $C_1(\varepsilon_0)$ is uniformly integrable.
				The first part of that theorem therefore gives the third statement.
			\end{proof}

\end{document}